\documentclass{m2an}
\usepackage[T1]{fontenc}
\usepackage{amsmath}
\usepackage{amssymb}
\usepackage{stmaryrd}
\usepackage{tikz}
\usepackage{pgfplots}
\usepackage{xcolor}
\usepackage{booktabs}
\usepackage{hyperref}
\usepackage{algorithmic,algorithm}
\usepackage[utf8]{inputenc}

\usepackage{csquotes}
\usepackage[french,english]{babel}

\usetikzlibrary{calc}
\usetikzlibrary{arrows,decorations.markings}
\usetikzlibrary{pgfplots.groupplots}

\usepackage{siunitx}

\newcommand\flot{\phi} 
\newcommand\inertial{\mathcal{I}}
\newcommand\rotating{\mathcal{R}}
\newcommand\Ext{\text{ext}}
\newcommand\vect[1]{\underline{#1}}

\newcommand\der[2]{\dfrac{\mathrm{d} #1}{\mathrm{d} #2}}
\newcommand\derpart[2]{\dfrac{\partial #1}{\partial #2}}

\newcommand\R{\mathbb{R}}

\let\leq\leqslant
\newcommand\D{\mathrm{d}}
\newcommand\norm[2]{\left\|#1\right\|^{#2}}
\newcommand\prodscal[2]{\langle#1,\;#2\rangle}
\newcommand\sphere{\mathbf{S}}

\DeclareMathOperator*{\argmin}{arg\,min}

\newcommand\Ftir{\mathcal{S}}
\newcommand\cost{\mathcal{C}}
\newcommand\Prob[1]{\mathcal{P}_{#1}}
\newcommand\abs[1]{\left|#1\right|}
\newcommand\OCP{(OCP)}
\newcommand\Ham{\mathcal{H}}

\newcommand\Lya{\text{Lya}}
\newcommand\Halo{\text{Halo}}
\newcommand\Het{\text{Het}}
\newcommand\Nat{\text{nat}}
\newcommand\Mani{\mathcal{M}}
\newcommand\tot{\text{tot}}
\newcommand\multi{\text{multi}}
\newcommand\thrust{\text{thrust}}

\newcommand\Tmax{T_{\text{max}}}
\newcommand\Isp{I_\text{sp}}
\newcommand\G{g_0}
\newcommand\latin[1]{\textit{#1}}

\newcommand\obj{\text{obj}}
\newcommand\energy{\mathcal{E}}

\newenvironment{Remark}{%
\medskip\par
\noindent \textbf{Remark:} }{\par\medskip}

\hypersetup{pageanchor=true}

\newcommand\pardesc[1]{\par\medskip\noindent\textbf{#1}~}

\begin{document}

\title{Low-Thrust Lyapunov to Lyapunov and Halo to Halo
  Missions with $L^{2}$-Minimization}\thanks{This work is dedicated to Philippe Augros.}

\author{Maxime Chupin}\address{
  LJLL-UPMC --  Paris, Airbus Defence and Space -- Les Mureaux, \email{chupin@ljll.math.upmc.fr}}

\author{Thomas Haberkorn}\address{MAPMO -- Orléans,
  \email{thomas.haberkorn@univ-orleans.fr}}

\author{Emmanuel Trélat}\address{LJLL-UPMC -- Paris,
  \email{emmanuel.trelat@upmc.fr}}

\date{Received: date / Accepted: date}

\begin{abstract}
In this work, we develop a new method to design energy minimum
low-thrust missions ($L^{2}$-minimization). In the Circular
Restricted Three Body Problem,
the knowledge of invariant manifolds helps us  initialize an indirect
method solving a transfer mission between periodic Lyapunov
orbits. Indeed, using the PMP,
the optimal control problem is solved using Newton-like algorithms
finding the zero of a shooting function. To compute a Lyapunov to
Lyapunov mission, we first compute an
admissible trajectory using a heteroclinic orbit between the two
periodic orbits. It is then used to initialize a multiple shooting method
in order to release the constraint. We finally optimize the terminal
points on the periodic orbits.  Moreover, we use continuation
methods on position and on thrust, in order to gain robustness. A more
general Halo to Halo mission, with different energies, is
computed in the last section
without heteroclinic
orbits but using invariant manifolds to initialize shooting
methods with a similar approach.
\end{abstract}

\begin{resume}
Dans ce travail, on développe une nouvelle méthode pour construire des
missions à faible poussée avec minimisation de la norme $L^{2}$ du
contrôle. Dans le problème circulaire restreint des trois corps,
la connaissance des variétés invariantes nous permet d'initialiser une
méthode indirecte utilisée pour calculer un transfert entre orbites
périodiques de Lyapunov. En effet, par l'application du Principe du
Maximum de Pontryagin, on obtient la commande optimale par le calcul
du zéro d'une fonction de tir, trouvé par un algorithme de Newton. Pour
construire la mission Lyapunov vers Lyapunov, dans
un premier temps, on calcule une trajectoire admissible en passant
par une trajectoire hétérocline reliant les deux orbites
périodiques. Celle-ci est alors utilisée pour initialiser un tir
multiple nous permettant de relacher la contrainte de rejoindre la
trajectoire hétérocline. Enfin, on optimise la position des points de
départ et d'arrivée sur les orbites périodiques. De plus, pour rendre
nos méthodes plus robustes, on utilise des méthodes de continuation
sur la position et sur la poussée. Dans la dernière section, on
contruit une mission plus générale Halo vers Halo avec des énergies
différentes. Cette fois, nous ne pouvons plus utiliser d'orbites
hétéroclines, mais on initialise la méthode de tir avec des
trajectoires des variétés invariantes de la même façon qu'avec
l'orbite hétérocline pour la mission Lyapunov vers Lyapunov.
\end{resume}



\maketitle


\section{Introduction}

Since the late '70s, study of libration point orbits has been of great
interest. Indeed,
several
missions such as ISEE-3 (NASA) in 1978, SOHO (ESA-NASA) in 1996,
GENESIS (NASA) in 2001, PLANK (ESA) in 2007 etc. have put this design
knowledge into practice. A more profound
understanding
of the available mission options has also emerged due to the theoretical,
analytical, and numerical
advances in many aspects of libration point mission design.

There exist a huge number of references on the problem of
determining low-cost trajectories by using the properties of Lagrange
equilibrium points.
For instance, the authors in \cite{Farquhar80,koon2000,Gomez00somezero,
  Zazzera2004} have
developed very efficient methods to find ``zero cost'' trajectories
between libration point orbits. Dynamical system methods are used
to construct heteroclinic orbits from invariant manifolds between
libration point orbits and
it allows to get infinite time uncontrolled
transfers. These orbits have been used with impulse engines of
spacecrafts to construct finite time transfers.
In this work, we want to perform the transfer with a low-thrust
propulsion,
so impulses to reach heteroclinic  orbits (or trajectories on
invariant manifolds) are prohibited.

Invariant manifolds have been used in a low-thrust mission
in~\cite{Mingotti2009, Mingotti2011}. The low-thrust propulsion is
introduced by
means of special attainable sets that are used in conjunction with
invariant manifolds to define a first-guess solution. Then, the
solution is optimized using an optimal control formalism.
One can note that~\cite{Mingotti2007} is the first work that combines
invariant manifolds and low-thrust in the
Earth-Moon system.

Much efforts have been dedicated to  the design of  efficient methods to
reach periodic orbits, Halo orbits, around equilibrium points in the
three body problem. For example, in~\cite{Senent2005, Ozimek2010},
authors use indirect methods and direct multiple shooting methods to reach an
insertion point on a manifold to reach asymptotically a Halo orbit
in the Earth-Moon system. Moreover, using transversality conditions,
 the position of the insertion point on the manifold is optimized on
 the manifold. Low-thrust,
stable-manifold
transfers to Halo orbits are also shown in~\cite{Starchville1997}.

On the same topic,
in~\cite{Mingotti2007}, authors use direct
methods to reach a point on a stable manifold of a Halo orbit from a
GTO orbit. A transfer from the Halo orbit to a
Lunar-Orbit is established as well. The $L_{2}$-norm of the control is
minimized by
a direct transcription and non linear programming.
In \cite{Martin2010},  the position of insertion point on the
manifold is optimized.

  We can notice that in the  interesting work~\cite{Daoud2011},
  indirect method combined with continuation methods have been used to
  design missions from an Earth Geostationary Orbit to a Lunar
  Orbit. Indeed,  continuations are used from the two body problem to
  the three body problem, the minimum time problem is studied and
  solved, and continuations
  between energy minimization and fuel consumption minimization are
  computed.

  Moreover, in~\cite{Francesco2015}, the developed methods  involve
  the minimum-time problem, the minimum energy problem and the minimum fuel
  problem to reach a fixed point on a Halo orbit starting from a
  periodic orbit around Earth. Continuations
  on the thrust are used as well as Newton and bisection methods
  (indirect methods).
  In these last two contributions, manifolds are not used to help solve the
  formulated problem.

In \cite{Epenoy}, the author recently developed an efficient
method to compute an optimal low-thrust transfer trajectory in
\emph{finite} time
without using invariant manifolds of the three
body problem. It is based on a three-step
solution method using indirect methods and continuations methods and it
gives good results.

  The philosophy of the method developed in this work is to use the
  natural dynamics as invariant manifolds, providing free parts for
  transfer, and to initialize a global multiple shooting method freeing
  the constraints to stay on the manifold. The invariant manifolds are
  just there to help obtain convergence for a shooting method.
  For references on techniques used in our work such as continuation on
cost,
smoothing techniques, optimization techniques one may
read~\cite{Haberkorn2004,Jiang2012,Bertrand2002}.

To compute the required transfer with a low thrust and minimizing the
energy ($L^{2}$-minimization), we will use
indirect methods coming from the Pontryagin Maximum
Principle~\cite{Pontryagin}.
Initialization of indirect methods with dynamics properties
is a real challenge in order to improve the efficiency of indirect
methods (see~\cite{Trelat2012} and references therein). Indeed, the
main difficulty of such
methods is to
initialize the Newton-like algorithm, and the understanding of the
dynamics can be very useful to construct an admissible trajectory for
the
initialization. Moreover, continuation methods as used
in~\cite{Haberkorn2006} or~\cite{CoCaGe2012}, are crucial to give robustness to these
indirect methods.
 The aim of this paper is to combine all these
mathematical aspects of dynamics, optimal control and continuation
methods to design low-thrust transfers
between libration point orbits.
  Principle, we only get necessary conditions of optimality. It
  would be interesting to check the second order sufficient
  conditions, with focal point tests for example.

The outline for the article is as follows. First, in
Section~\ref{sec:mission}, we introduce the mission we want to
perform and compute, introduce the paradigm of the
Circular Restricted Three Body Problem and  state our optimal
control problem.

Then, in Section~\ref{sec:dynamics}, we  recall dynamical
properties of the circular restricted
three body problem, such as equilibrium points, Lyapunov periodic
orbits, and invariant manifolds. We present the mathematical tools
used to numerically compute the periodic orbits and the manifolds.  In
particular, we  introduce in this part
the continuation method that we will use throughout this article.

Then, in Section~\ref{sec:method}, we  develop our method with an
example  mission. We first
compute a heteroclinic orbit between the two Lyapunov periodic
orbits. Then, fixing the departure point near $L_{1}$ and the arrival
point near $L_{2}$ and  with a not to small thrust (\SI{60}{N}), we
perform two small
transfers from the Lyapunov
orbit around $L_{1}$ to the heteroclinic one, and from the
heteroclinic orbit to
the Lyapunov orbit around $L_{2}$. Then, thanks to a multiple shooting
method  we release the constraint on the position of the matching
points  on the heteroclinic orbit and
decrease
the thrust to the targeted one (\SI{0.3}{N}). Finally, we optimize the
 departure and  arrival points on the periodic orbits to satisfy the
necessary transversality conditions given by the Pontryagin Maximum
Principle. In section~\ref{sec:mission2}, we present another mission with
a heteroclinic
orbit with two revolutions around the Moon.

Finally, in section~\ref{sec:HaloHalo}, we apply the method
  to a more general mission: a Halo to Halo mission for two
  periodic orbits with different energies. In this case, there is no heteroclinic
  orbit, so we construct an admissible trajectory with
  5 parts. Two of them are trajectories on invariant manifolds, and
  the three others are local transfers: 1/ from one of the  Halo orbits to a free
  trajectory, 2/ between both free trajectories and 3/ from
  the second free trajectory to the second Halo orbit. Thanks to this
  five part admissible trajectory, we are able to initialize a multiple
  shooting method that computes  an optimal trajectory (which is not
  constrained to reach any invariant manifolds). As previously, we
  optimize the terminal points on Halo orbits.

\section{The Mission}\label{sec:mission}

\subsection{Circular Restricted Three Body Problem (CRTBP)}

We use the paradigm
of the Circular Restricted Three Body Problem. In this section we will
follow the description by \cite{Koon2006}.

Let us consider a spacecraft in the field of attraction of  Earth
and Moon.  We consider an inertial frame $\inertial$ in which the vector
differential equation for the spacecraft’s motion is written as:
\begin{equation}\label{eq:3BPmotionInertial}
  m\der{\vect{R}}{t} = -GM_1 m \dfrac{\vect{R}_{13}}{R^3_{13}}-GM_2 m \dfrac{\vect{R}_{23}}{R^3_{23}}
\end{equation}
where $M_1$, $M_2$ and $m$ are the masses respectively of Earth,
the Moon and the spacecraft, $\vect{R}$ is the spacecraft vector
position, $\vect{R}_{13}$ is the vector Earth-spacecraft and
$\vect{R}_{23}$ is the vector Moon-spacecraft. $G$ is the gravitational
constant. Let us describe the simplified general framework we will use.

\pardesc{Problem Description.} To simplify the problem, and use a
general framework, we consider
the motion of the spacecraft $P$ of negligible  mass  moving under the
gravitational influence of the two
masses $M_{1}$ and $M_{2}$, referred to as the primary masses, or
simply the \emph{primaries} (here Earth and  Moon). We denote
these primaries by $P_{1}$ and
$P_{2}$. We assume that the primaries have circular orbits around their
common center of mass. The particle $P$ is free to move all around
the primaries but cannot affect their motion.

The system is made adimensional by the following choice of units:
the unit of mass is taken to be $M_{1}+M_{2}$; the unit of length is
chosen to be the constant distance between $P_{1}$ and $P_{2}$; the
unit of time is chosen such that the orbital period of $P_{1}$ and
$P_{2}$ about their center of mass is $2\pi$. The universal constant
of gravitation then becomes $G=1$. Conversions from units of
distance, velocity and time in the unprimed, normalized system to the
primed, dimensionalized system are
\begin{equation}\label{eq:normalized}
    \text{distance}\quad  d'=l_*d,\quad
    \text{velocity} \quad s'=v_*s,\quad
    \text{time} \quad t'=\frac{t_*}{2\pi}t,
\end{equation}
where we denote by $l_*$ the distance between $P_{1}$ and $P_{2}$, $v_*$ the
orbital velocity of $P_{1}$ and $t_*$ the orbital period of $P_{1}$ and
$P_{2}$.

We define the only parameter of this system as
\begin{equation*}
  \label{eq:mudef}
  \mu=\dfrac{M_{2}}{M_{1}+M_{2}},
\end{equation*}
and call it the \emph{mass parameter}, assuming that $M_{1}>M_{2}$.

In table \ref{tab:parameters}, we summarize the  values of
all the constants for the Earth-Moon CR3PB for numerical computations.

\begin{table}[ht]
  \centering
  \begin{tabular}{ccccc}
    \toprule
    System & $\mu$ & $l_*$ & $v_*$ & $t_*$ \\
    \midrule
    Earth-Moon & \num{1.215e-2} & \SI{384402e03}{\kilo\meter}&
    \SI{1.025}{\kilo\meter\per\second} & \SI{2.361e06}{\second}\\
    \bottomrule
  \end{tabular}
  \caption{Table of the parameter values for the Earth-Moon system.}
  \label{tab:parameters}
\end{table}

\pardesc{Equations of Motion.} If we write the equations of motion in a
rotating frame $\rotating$ in which the two primaries are fixed (the
angular velocity is the angular velocity of their rotation around
their center of mass, see~\cite{Daoud2011}), we obtain that the
coordinates of
$P_{1}$ and $P_{2}$ are respectively
$\chi_{P_{1}}=(-\mu,0,0,0,0),$ and $\chi_{P_{2}}=(1-\mu,0,0,0,0)$.
Let us call $x_{1}^{0}=-\mu$ and
$x_{2}^{0}=1-\mu$, and by writing the state $\chi=\left(x,y,\dot{x},\dot{y}\right)^T=\left(x_1,x_2,x_3,x_4\right)^T,$
we obtain
\begin{equation}
  \label{eq:motion}
  \left\{%
    \begin{aligned}
      \dot{x}_1&=x_3\\
      \dot{x}_2&=x_4\\
      \dot{x}_4&=x_1+2x_4-(1-\mu)\frac{x_1-x_1^0}{r_1^3}-\mu\frac{x_1-x_2^0}{r_2^3}\\
      \dot{x}_5&=x_2-2x_3-(1-\mu)\frac{x_2}{r_1^3}-\mu\frac{x_2}{r_2^3}\\
    \end{aligned}
  \right.
\end{equation}
where
$r_1=\sqrt{\left(x_1-x_1^0\right)^2+x_2^2}$ and $r_2=\sqrt{\left(x_1-x_2^0\right)^2+x_2^2}$
are respectively the distances between $P$ and primaries $P_1$ and $P_2$.

We can define the potential
$U(x_1,x_2)=-\frac{1}{2}\left(x_1^2+x_2^2\right)-\frac{1-\mu}{r_1}-\frac{\mu}{r_2}-\frac{1}{2}\mu\left(1-\mu\right).$
We denote by $F_{0}$  the vector field of the
system and we define the energy of a state point as
\begin{equation}
  \label{eq:Energy}
  \energy(\chi) = \frac{1}{2}(\dot x^{2}+\dot y^{2})+U(x,y).
\end{equation}
Note that the energy is constant as the system evolves over time
(conservation law).

In the first part of this work, we will consider a planar motion,
hence, we only have a
$\R^{4}$-state in the orbital plane of the primaries, but this can be
easily extended to the spatial case.

\subsection{Design of the Transfer}

We want to design a mission going from a periodic Lyapunov orbit
around $L_1$
to a periodic  Lyapunov orbit around $L_2$ using a low-thrust engine
in the
Earth-Moon system (see figure~\ref{fig:U2U3}). A full
description
of these periodic orbits is given in section~\ref{sec:Lyapunov}.  In
order to
perform such a mission, we will use the properties
introduced in section~\ref{sec:manifolds}: the invariant
manifolds. Indeed, if we are able to find
an intersection between an ``$L_1$ unstable manifold'' and an ``$L_2$
stable manifold'', we get an asymptotic trajectory that performs the
mission with a zero thrust, called a heteroclinic orbit (see
section~\ref{sec:het}).

In the classical literature, such a mission is usually designed by using impulse
to reach the heteroclinic orbit from the Lyapunov orbit around $L_{1}$
and then another impulse to reach the Lyapunov orbit from the
heteroclinic one. Since we design a low-thrust transfer, following
this method is unrealistic. In~\cite{Epenoy}, the author developed a
three-step method to perform a low thrust low energy
trajectory between Lyapunov orbits of the same energy without using
invariant manifolds. At his first step, he uses a feasible
quadratic-zero-quadratic control structure  to initialize his
method. In this work we will use the knowledge of a zero cost
trajectory, the heteroclinic orbit, to initialize an indirect shooting
method (Newton-like method for optimal control problem) provided by
applying the Pontryagin Maximum Principle.

\subsection{Controlled Dynamics}

We first describe the model for the evolution of our spacecraft in the
CRTBP. In non normalized coordinates (see~\eqref{eq:3BPmotionInertial}), the
controlled dynamics is the $m\der{\vect{R}}{t} = -GM_1 m \dfrac{\vect{R}_{13}}{R^3_{13}}-GM_2 m \dfrac{\vect{R}_{23}}{R^3_{23}}+T(t),$
where $T$ is the spacecraft driving force, and $m$ is the time
dependant
mass of the spacecraft. The equation for the evolution of the mass
is
\begin{equation*}
  \label{eq:massevol}
  \dot m(t) = -\beta \norm{T(t)}{},
\end{equation*}
where $\beta$ is computed with the two parameters
$\Isp$ and $\G$. Specific impulse ($\Isp$) is a
measure of the efficiency of rocket and jet engines. $\G$ is the
acceleration at Earth's surface. The inverse of
the average exhaust speed, $\beta$, is equal to $\frac{1}{\Isp\G}$.
Moreover, the thrust is constrained by  $\norm{T(t)}{}\leq \Tmax$  for
all $t$.

Using the normalization parameters~\eqref{eq:normalized}, denoting by
$\beta_{*}$ the normalized parameter $\beta$ initially in
\si{m^{-1}.s}, the mass evolution is
\begin{equation*}
  \label{eq:massevol2}
      \dot m(t) = -\beta_{*} \frac{t_*^2}{4\pi^2 l_*}\Tmax \norm{u(t)}{}
 \end{equation*}
Moreover, we introduce the control $u$ such that $\abs{u(t)}\leq 1$
and denote the normalized coefficient
$\frac{t_*^2}{4\pi^2 l_*}\Tmax$ by $\epsilon$.

In short, we write the system as :
\begin{equation*}
  \label{eq:vectorfield}
  \left\{
    \begin{array}{l}\displaystyle
      \dot x = F_{0}(x)+\frac{\epsilon}{m}\sum_{i=1}^{2}u_{i}F_{i}(x),\\
      \displaystyle\dot m = -\beta_{*}\epsilon\norm{u}{},
    \end{array}
  \right.
\end{equation*}
where $F_{0}$ is the natural vector field defined
by~\eqref{eq:motion}. Here  we have
$
F_{1}(x)=\begin{pmatrix}0\\0\\1\\0
\end{pmatrix}$, and $F_{2}(x)=\begin{pmatrix}0\\0\\0\\1
\end{pmatrix}$.
This can be easily extended to the spatial case.

\pardesc{Controllability.} In \cite{Caillau2012}, it is proved that the
CRTBP with a non evolving mass is controllable for a suitable
subregion of the phase-space, denoted by $X_{\mu}^{1}$, where the
energy is greater than the energy of $L_{1}$:
\begin{thrm}
  For any $\mu\in(0, 1)$, for any positive $\epsilon$, the circular restricted
  three-body problem is controllable on $X_{\mu}^{1}$.
\end{thrm}
Using proposition (2.2) in~\cite{CaillauThese}, one can extend
this result to the system with an evolving mass.

\subsection{Optimal Control Problem \OCP}

Our main goal in this work is to solve an optimal control problem. We want
to go from the Lyapunov orbit around $L_{1}$ to the Lyapunov orbit
around $L_{2}$ with minimal energy. Mathematically we write this
problem as follows
\begin{equation}
  \label{eq:OCPgoal}
  \Prob{g}\left\{
    \begin{array}{l}
      \displaystyle \cost_{g} = \min \int_{0}^{t_{f}}\norm{u}{2}\D t,\\
      \displaystyle\dot x = F_{0}(x)+\frac{\epsilon}{m}\sum_{i=1}^{2}u_{i}F_{i}(x),\\
      \displaystyle\dot m = -\beta_{*}\epsilon\norm{u}{},\\
      \norm{u}{}\leq 1,\\
      x(0)\in\Lya_{1},\text{ and }x(t_{f})\in\Lya_{2}.
    \end{array}
  \right.
\end{equation}

Let us summarize  the steps in the method we developed to solve this
problem :
\begin{enumerate}
\item First, we find a heteroclinic orbit from the Lyapunov orbit
  around $L_{1}$ to the Lyapunov orbit around $L_{2}$.
\item Then, we realize a short transfer from \emph{a fixed point} on
  the Lyapunov orbit
  around $L_{1}$ to the heteroclinic orbit.
\item Similarly, we realize a transfer from the heteroclinic orbit
  to \emph{a fixed point} on the Lyapunov orbit around $L_{2}$.
\item Then we release the constraint on the position of the matching
  connections on the heteroclinic
  orbit using a multiple shooting method and we decrease the maximal
  thrust.
\item Finally, we optimize the position of the two fixed points on
  $\Lya_{1}$ and $\Lya_{2}$ to satisfy the transversality condition for
  problem~\eqref{eq:OCPgoal}.
\end{enumerate}

We note that in steps 2 to 4  (where we are solving optimal
control problems), we have fixed the departure and arrival points to
simplify the problem. The last step consists in releasing these
constraints.

  \begin{Remark}The \emph{real} problem that we want to solve is
  the minimization of the consumption of fuel (the maximization of
  the final mass). This is done by
  considering the minimization of the $L_{1}$-norm of $u$
  \[
  \cost_{g}^{L_{1}} = \min \int_{0}^{t_{f}}\norm{u}{}\D t.
  \]
  Unfortunately, this implies numerical difficulties and for
  simplicity, we only consider here the $L_{2}$-minimization
  problem. One can see \cite{Daoud2011}, or \cite{Francesco2015} where
  the authors attempt to consider the $L^{1}$-minimization.
  This is one of the perspective of this work using for example
  another continuation on the cost.
\end{Remark}

\section{Properties of CRTBP}\label{sec:dynamics}

In this Section, we recall some properties of the CRTBP. In
particular, we introduce equilibrium points, Lyapunov orbits and
invariant manifolds. We explain how to numerically compute these
orbits (see Section~\ref{sec:family}). We have improved the method
used in~\cite{Archambeau2011} using the energy as continuation
parameter. Finally, we introduce the invariant manifolds and how we
can get a numerical approximation.

\subsection{Lyapunov Orbits}\label{sec:Lyapunov}

\begin{figure}[ht]
  \centering
  \includegraphics[page=2]{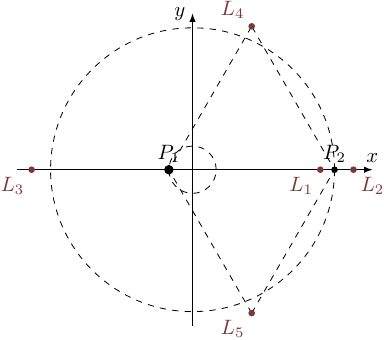}
  \caption{Localization of Lagrange's points.}
  \label{fig:Lagrange}
\end{figure}

\pardesc{Equilibrium Points.} The Lagrange points are the equilibrium
points of the circular restricted three-body problem.
Euler~\cite{Euler} and Lagrange~\cite{Lagrange} proved the existence of five equilibrium
points: three collinear
points on the axis joining the center of the two primaries, generally
denoted by $L_{1}$, $L_{2}$ and $L_{3}$, and
two equilateral points denoted by $L_{4}$ and $L_{5}$ (see figure~\ref{fig:Lagrange}).

Computing equilateral points $L_{4}$ and $L_{5}$ is not very
complicated, but it is not possible to find exact solutions for
collinear equilibria $L_{1}$, $L_{2}$ and $L_{3}$.  We refer
to~\cite{Szebehely}, for series expressions. We recall
that  the collinear points are shown to be unstable (in every system),
whereas $L_{4}$ and $L_{5}$ are proved to be stable under some
conditions (see \cite{meyer2010}).

\pardesc{Periodic Orbits.}
The Lyapunov center theorem ensures the existence of
periodic orbits around equilibrium points (see \cite{meyer2010,
  Bonnard2009} and references therein). To use this theorem, one has
to linearize the system and compute the eigenvalues of the linearized system.

 In the planar case, applying this
theorem to the collinear points $L_{1}$, $L_{2}$ and $L_{3}$ we get a
one-parameter family of periodic orbits around each one. These
periodic orbits are called Lyapunov orbits
and are homeomorphic to a circle. In this work, we denote by $\Lya_{i}$ a
Lyapunov orbit around the equilibrium point $L_{i}$. For the spatial
case, periodic orbits are
called Halo orbits or Lissajous orbits
(see e.g.~\cite{Gomez1997}).

\pardesc{Numerical computation.} We will describe the method to
compute Lyapunov orbits around collinear Lagrange points. For the
spatial case we follow the same method. To
find these periodic orbits, we use a Newton-like method. Since
equations in the coordinate system centered on $L_{i}$ are symmetric, if
we
consider a periodic solution
$\chi(t)=(x(t),y(t),\dot x(t),\dot y(t))$ of period $t_{\chi}$, then there
exists $t_{0}$ such that
\[
\left\{
  \begin{array}{l}
    x(t_{0}) = x_{0},\\
    y(t_{0}) = 0,\\
    \dot x(t_{0}) = 0,\\
    \dot y(t_{0}) = \dot y_{0},
  \end{array}
\right.
\text{ and }
\left\{
  \begin{array}{l}
    x(t_{0}+t_{\chi}/2) = x_{1},\\
    y(t_{0}+t_{\chi}/2) = 0,\\
    \dot x(t_{0}+t_{\chi}/2) = 0,\\
    \dot y(t_{0}+t_{\chi}/2) = \dot y_{1}.
  \end{array}
\right.
\]

Since $t_{0}$ could be chosen to be equal to zero and fixing
$x_{0}$, we just
have to find $(\dot y_{0}, t_{\chi})$ such that, denoting by $\flot$ the
flow of the dynamical system, and  $\chi_{0}=(x_{0},0,0,\dot y_{0})$,
the function $\Ftir_{L}$ satisfies:
\begin{equation}\label{eq:Ftir}
  \Ftir_{L}(t_{\chi},\dot y_{0})=\begin{pmatrix}\flot_{2}(t_{\chi}/2,\chi_{0})\\\flot_{3}(t_{\chi}/2,\chi_{0})
  \end{pmatrix}
  =\begin{pmatrix}0\\0
  \end{pmatrix}.
\end{equation}

In practice, we fix for example the value of $x_{0}$ (respectively of  $z_{0}$
in the spatial case) in order to be
left with  finding a
zero of a function of two variables $(\dot y_{0},
t_{\chi})$ in $\R^{2}$ . Obviously, we can extend this to a periodic
orbit in
$\R^{6}$.

The main difficulty is to initialize the Newton-like algorithm. The idea
is to find an analytical approximation of the orbit to a certain
order, and then inject this into the Newton-like algorithm. In this work, and
because the Lyapunov orbit is not very difficult to compute, we follow
\cite{Richardson}. For various orbits in $\R^{6}$, see
\cite{Archambeau2011,Jorba1997,Farquhar73} and references therein.

\subsection{Computing the family}\label{sec:family}

In order to use these orbits to construct
the targeted  mission, it is very useful to be able to compute the family
of periodic orbits, providing us with different orbits that have
different energies.

\subsubsection{Continuation methods}\label{sec:continuation}

To explain how we get the family of periodic orbits, let us
introduce continuation methods, for a more complete introduction,
see~\cite{Allgower1990}. The main idea is to construct a
family of problems denoted by
$\left(\Prob{\lambda}\right)_{\lambda\in[0,1]}$ indexed by a
parameter $\lambda\in[0,1]$. The initial problem $\Prob{0}$ is
supposed to be easy to solve, and the final problem $\Prob{1}$ is the
one we want  to solve.

Let us assume that we have solved numerically $\Prob{0}$, and consider a
subdivision $0=\lambda_{0}<\lambda_{1}<\cdots<\lambda_{p}=1$ of the
interval $[0,1]$. The solution of $\Prob{0}$ can be used to initialize
the Newton-like method applied to $\Prob{\lambda_{1}}$. And so on, step by
step, we use the solution of $\Prob{\lambda_{i-1}}$ to initialize
$\Prob{\lambda_{i}}$. Of course, the sequence
$\left(\lambda_{i}\right)$ has to be well chosen and eventually should
be refined.

Mathematically, for this method to converge, we need that
the family of problems to depend continuously on the parameter
$\lambda$.
See~\cite[chap. 9]{Bonnard2009} for some justification of the
method.

From the numerical point of view, there exist many methods and
strategies for implementing continuation or homotopy methods.
We can distinguish between differential pathfollowing, simplicial
methods, predictor-corrector methods, etc.
In this work, we implement a predictor-corrector method because it
is suitable for our problem. Here, we use a ``constant''
prediction: the solution of problem $\Prob{\lambda_{i-1}}$ is used to initialized
the resolution of problem $\Prob{\lambda_{i}}$. We can
note that there exist many codes which can be found on the web, such
as the well-known Hompack90~\cite{hompack90} or
Hampath~\cite{CoCaGe2012}. For a survey about different
results, challenges and issues on continuation methods,
see~\cite{Trelat2012}.

\subsubsection{Application to the family of orbits}\label{sec:familyorbit}

Since we had to choose a parameter $x_{0}$ to write the zero function
$\Ftir_{L}$ in \eqref{eq:Ftir}, it is natural to use this
parameter to perform the continuation that computes the family of
orbits. Indeed, we can choose to reach a certain
$x_{0}^{\obj}$ (respectively a so called  \emph{excursion} $z_{0}^{\obj}$ in
the spatial case). So
we can define our continuation as :
\[
\Prob{\lambda}:\left\{\begin{array}{l}
  \text{for } x_{0}^{\lambda}=(1-\lambda)x_{0}+\lambda
  x_{0}^{\obj}\\
  \Ftir_{L}^{\lambda}(t_{\chi},\dot y_{0})=\begin{pmatrix}\flot_{2}(t_{\chi}/2,\chi_{0}^{\lambda})\\\flot_{3}(t_{\chi}/2,\chi_{0}^{\lambda})
  \end{pmatrix}
  =\begin{pmatrix}0\\0
  \end{pmatrix}
\end{array}
\right.
\]
where $\chi_{0}^{\lambda}=(x_{0}^{\lambda},0,0,\dot y_{0})$. Thanks to
the analytical approximation provided by \cite{Richardson} or
\cite{Jorba1997}, we can solve the initial problem
$\Prob{0}$. We can note that such analytical approximation does not
work for every $x_{0}$. Using the continuation method described
previously, we can
get a family of periodic orbit.

However, for some periodic orbits (Halo family), we can observe that
the continuation fails when we converge to the equilibrium
point ($x_{0}^{\obj}\to 0$). A much better continuation parameter is
\emph{energy}. It releases the constraint on the
parameter $x_{0}$, and allows us to reach any periodic
orbit, in particular the algorithm converges to the energy of
$L_{i}$, $i\in\{1,\dots,3\}$. Moreover, it is a significantly more natural
parameter, keeping in
mind the fact that we will construct a controlled transfer
method. Section~\ref{sec:manifolds} will provide an extra argument in
favor of the
energy parameter. It seems to be the first time that this
  continuation is done with the energy as the continuation
  parameter. This avoids numerical problems when reaching energy close
  to $L_{i}$.

Thanks to the analytical approximation, we get a first periodic
orbit with  energy $\energy_{0}$, and we want to reach a
prescribed energy $\energy_{1}$ so we define the following
family of problems:
\begin{equation*}\label{eq:contOrbite}
\Prob{\lambda}^{\energy}:\begin{array}{l}
    \Ftir_{\energy}^{\lambda}(t_{\chi}, x_{0},\dot
    y_{0})=\begin{pmatrix}\flot_{2}(t_{\chi}/2,\chi_{0})\\\flot_{3}(t_{\chi}/2,\chi_{0})\\ \energy(\chi_{0})-\energy_{\lambda}
    \end{pmatrix}
    =\begin{pmatrix}0\\0\\0
    \end{pmatrix}
  \end{array}
\end{equation*}
where $\energy(\chi_{0})$ is the energy of the trajectory starting at
$\chi_{0}$ and
$
\energy_{\lambda}=(1-\lambda)\energy_{0}+\lambda\energy_{1}$.

In the continuation, we just use a
predictor-corrector continuation with a ``constant'' prediction as
explain before
(see~\cite{Allgower1990}). If necessary, one could improve this by using a
  linear predictor continuation, but using energy as a continuation
  parameter, continuation was very fast and easy, and did not
  require improvements.

\begin{figure}[ht]\centering
  \includegraphics[page=3]{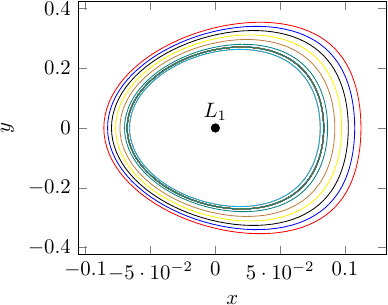}
  \caption{A family of Lyapunov orbits around $L_{1}$ in the
    Earth-Moon system (Richardson coordinates)}
\label{fig:familyLya}
\end{figure}

Figure~\ref{fig:familyLya} shows an example of a family of
Lyapunov orbits around $L_{1}$ in the Earth-Moon system.

\subsection{Invariant Manifolds}\label{sec:manifolds}

All the periodic orbits described in the previous section come with
their invariant manifolds, that is to say, the sets of phase points
from which the trajectory converges to the periodic orbit, forward for
the \emph{stable} manifold and backward for the \emph{unstable}
manifold. These manifolds can be very useful to design interplanetary
missions because as separatrix, they are some sort of gravitational
currents. We refer to \cite[chap. 4]{Koon2006} for the proof of
existence and
a more detailed explanation of these manifolds.
For the sake of numerical reproducibility, we recall some  well
known properties.

\pardesc{Monodromy Matrix.} We introduce a tool of dynamical
systems: the monodromy matrix. Some properties of this matrix
are needed to numerically compute the invariant manifolds. For more
details, see~\cite{meyer2010, Koon2006}.

Let $\bar{x}(\cdot)$ be a periodic solution of the dynamical system
with  period $T$ and $\bar{x}(0)=\bar x_{0}$. Denoting by $\flot$ the
flow of the system, the monodromy matrix
$M$ of the periodic orbit for the point $\bar x_{0}$ is defined as
\begin{equation*}\label{eq:monodromy}
  M=\derpart{\flot(T;\bar x_{0})}{x_{0}}.
\end{equation*}
It determines whether initial perturbations $\delta\bar x_{0}$ of
the periodic orbit decay or grow.

\pardesc{Local Approximation to Compute Invariant Manifolds.} Using
the Poincaré map we can show that the eigenvectors corresponding to
eigenvalues of the monodromy matrix are linear approximations of the
invariant manifolds of the periodic orbit. For the planar Lyapunov
orbits in the CRTBP, we show that the four eigenvalues of $M$ are
$
\lambda_{1}>1,\quad \lambda_{2}=\frac{1}{\lambda_{1}}, \quad \lambda_{3}=\lambda_{4}=1.
$
The eigenvector associated with eigenvalue $\lambda_{1}$ is in the
unstable direction and the eigenvector associated with  eigenvalue
$\lambda_{2}$ is in the stable direction.

\begin{figure}[ht]
  \includegraphics[page=4]{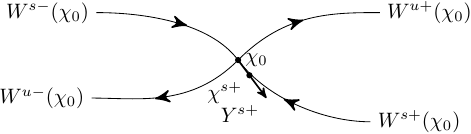}
  \caption{Illustration of the method to compute invariant manifolds.}
  \label{fig:manifoldX}
\end{figure}

Then, the method to compute invariant manifolds is the following:
\begin{enumerate}
\item First, for $\chi_{0}$ a point on the periodic orbit, we compute the
  monodromy matrix and its eigenvectors. Let us denote by $Y^{s}(\chi_{0})$
  the normalized stable eigenvector and by $Y^{u}(\chi_{0})$ the normalized
  unstable eigenvector.
\item Then,  let
  \begin{equation}
    \label{eq:pertubLya}
    \begin{array}{l}
      \chi^{s\pm}(\chi_{0})=\chi_{0}\pm\alpha Y^{s}(\chi_{0}),\\
      \chi^{u\pm}(\chi_{0})=\chi_{0}\pm\alpha Y^{u}(\chi_{0}),
    \end{array}
  \end{equation}
  be the initial guesses for (respectively) the stable and unstable
  manifolds. The magnitude of $\alpha$ should be small enough to be
  within the validity of  the linear estimate but not too small to
  keep a reasonable time of escape or convergence (for instance,
  see~\cite{Gomez1991a} for a discussion on the value of $\alpha$).
\item Finally, we integrate numerically the unstable vector forward
  in time, using both $\alpha$ and $-\alpha$ to generate the two
  branches of the unstable manifold denoted by $W^{u\pm}(\chi_{0})$. We do
  the same for the stable vector backwards, and we get the two
  branches of stable manifold $W^{s\pm}(\chi_{0})$ (see
  figure~\ref{fig:manifoldX}).
\end{enumerate}

\begin{figure}[ht]\centering
  \includegraphics{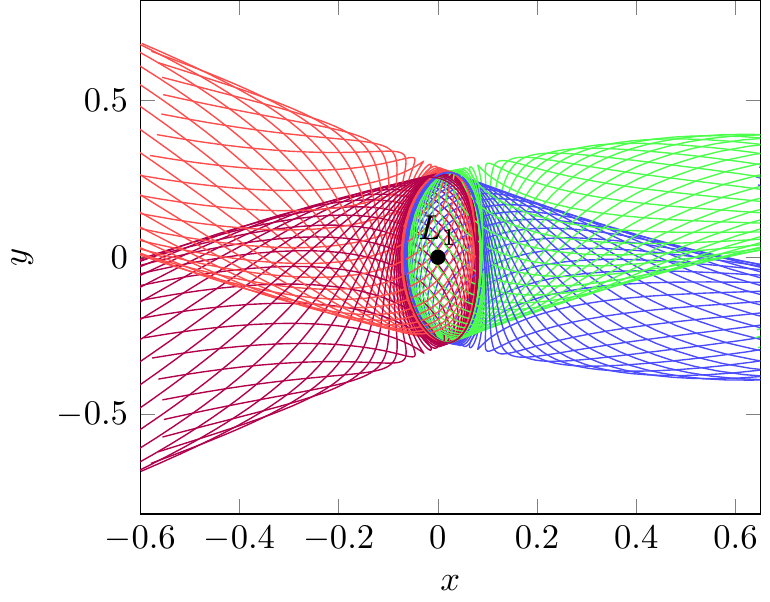}
  \caption{Manifolds of a Lyapunov orbit around $L_{1}$ of the
    Earth-Moon system and Richardson coordinates. The energy of these orbit and manifolds is
    $-1.59208$ in normalized coordinates (centered on the barycenter
    of the two primaries).}
  \label{fig:TubeL1}
\end{figure}

Following this process, we are able to compute the invariant manifolds
of any Lyapunov orbit at any energy (greater than the energy of
$L_{i}$, for an explanation of that, see the section about Hill regions
in~\cite{Koon2006}). We have represented parts of these manifold
for a Lyapunov orbit in the Earth-Moon system at  energy  $-1.59208$
in normalized coordinates in figure~\ref{fig:TubeL1}. Note that an
interesting study of fast numerical approximation of invariant
manifolds can be found in~\cite{Topputo201689}.

\begin{Remark} Since we are following invariant manifolds converging in
  infinite time to the periodic orbits (backward or forward), and
  because we are doing it numerically and so with a certain
  approximation, there exist long times for which we cannot obtain
  convergence. We have to tune the parameter $\alpha$
  in~\eqref{eq:pertubLya}  (we give in
  section~\ref{sec:het} the
  choice of the numerical value). Moreover, the multiple shooting method
  allows subdividing the time and keep each part to a reasonable time of
  integration.
\end{Remark}

\section{Constructing the Mission}\label{sec:method}

In this section we  explain all the steps of our method for solving
the problem~\eqref{eq:OCPgoal}. We first find a heteroclinic orbit
between the two Lyapunov orbits. Then we perform two short transfers
from $\Lya_{1}$ to the heteroclinic orbit, and from the heteroclinic
orbit to $\Lya_{2}$. Then, with a multiple shooting method we
release the constraint on the position of the matching connections on
the heteroclinic orbit. Finally, we
optimize the departure and arrival points previously fixed to simplify
the problem.

\subsection{The Heteroclinic Orbit}\label{sec:het}
Let us first find the heteroclinic orbit between a Lyapunov orbit around
$L_1$ and a Lyapunov orbit around $L_2$. One condition to be able to
find such an orbit is to compute an intersection between two
manifolds. Hence, these two manifolds should have the same energy. Since
the manifold and the Lyapunov orbit have the same energy, we must
compute two Lyapunov orbits around $L_1$ and $L_2$ with a given energy.

The study of the well known Hill regions (see \cite{Koon2006} and
references therein),
\latin{i.e.} the projection of the
energy surface of the uncontrolled dynamics onto the position space gives
us an indication of the interval of energy we can use. Indeed, we have to
compute an orbit with an energy greater than the $L_{2}$ energy. And
because we want to realize a low-thrust transfer, we choose to keep
a low energy. Moreover, we have a smaller region of possible
motion, and so, a possibly shorter transfer.

Using the method described in~\ref{sec:familyorbit}, we choose to get
two orbits with an energy of $-1.592081$ in the normalized
system.

\begin{figure}[ht]\centering
  \includegraphics[page=5]{pdffigures-crop}
  \includegraphics[width=0.5\textwidth]{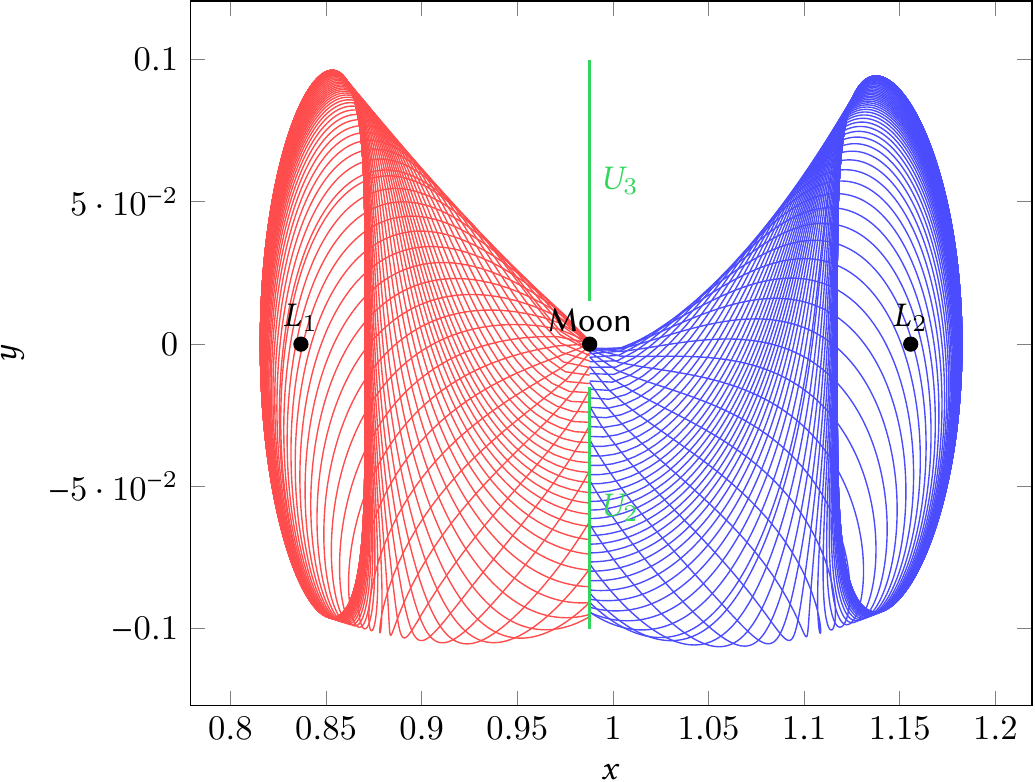}
  \caption{On the left: Planes $U_2$ and $U_3$ in the Earth-Moon
    system. On the right, unstable (red) and stable (blue) manifolds
    respectively from $L_{1}$ and $L_{2}$ stopping at the plane $U_{2}$}
  \label{fig:U2U3}
\end{figure}

\pardesc{Finding the intersection.}
To find an intersection, we introduce two 2D sections
$U_2=\{(x,y)\in \R^2,\; x=1-\mu,\; y<0\},$ and $U_3=\{(x,y)\in \R^2,\; x=1-\mu,\; y>0\}$.
We represent them in figure~\ref{fig:U2U3}.

Then, we compute the intersection of the unstable manifold from $L_1$ and
the stable manifold from $L_2$ with the space $U_2$ (of course, we
can do the
symmetric counterpart: stable manifold from $L_1$ and unstable manifold
from $L_2$ with the space $U_3$). Since the $x$-coordinate is
fixed by $U_2$ and because the energies of the two manifolds are equal, we
just have to compute the intersection in the $(y,\dot y)$-plan
(values of $\dot x$ are deduced from the energy
equation~\eqref{eq:Energy}).

\begin{figure}[ht]\centering
  \includegraphics[page=6]{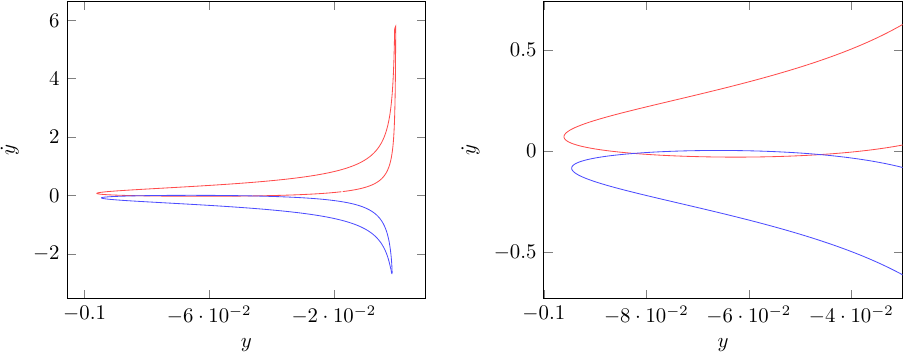}
  \caption{Section in the plane $U_{2}$ of a unstable manifold from
    $L1$ and a stable manifold from $L_{2}$. The energy is
    $-1.592081$. On the right, a zoom on the interesting area.}
  \label{fig:U2section}
\end{figure}

We show in figure~\ref{fig:U2section} the $U_{2}$-section and the
existence  of intersections for our particular energy. To find
precisely one intersection point, we have used once more a Newton-like
method. We can parametrize the section of one manifold with $U_{2}$
with only one parameter, the parameter of the Lyapunov orbit. We
denote by $\flot^{+}_{x=1-\mu}$, the flow propagating forward a state
point from $\Lya_{1}$ onto the space $U_{2}$, and by
$\flot^{-}_{x=1-\mu}$ the flow
propagating backward a state point from $\Lya_{2}$ onto the plan
$U_{2}$. Time of
propagation is fixed by the condition $x=1-\mu$.

We want to find  two points $\chi_{L_{1}}\in \Lya_{1}$ and
$\chi_{L_{2}}\in\Lya_{2}$
such that
\[\flot^{+}_{x=1-\mu}(\chi^{u+}(\chi_{L_{1}}))-\flot^{-}_{x=1-\mu}(\chi^{s-}(\chi_{L_{2}}))=0,\]
where $\chi^{u+}$ and $\chi^{s-}$ are defined in~\eqref{eq:pertubLya}.
This is an  equality in $\R^{2}$, and because each of the Lyapunov orbits is
parametrized with a one dimensional parameter (the time), our problem
is well posed.

To initialize the method we use a discretisation (100 points in this
particular example) of the Lyapunov orbits
and we take the two points minimizing the Euclidean norm in the $U_{2}$
section.

\begin{figure}[ht]\centering
  \includegraphics[page=7]{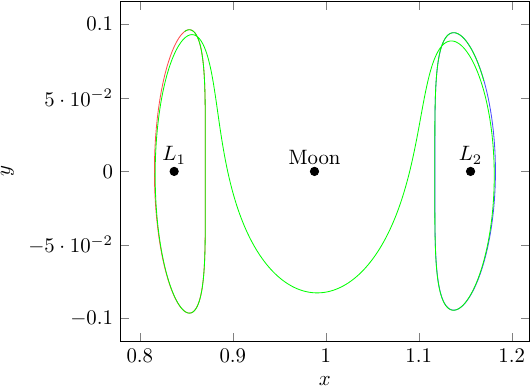}
  \caption{Heteroclinic orbit between two Lyapunov orbits in the
    Earth-Moon system. We get a travel time of $8.9613933501964$
    (normalized time) or $38.974$ days.}
  \label{fig:heteroclinic}
\end{figure}

In our case, with a value of energy equal to $-1.592081$ and
$\alpha=\frac{1}{384402}$ from~\eqref{eq:pertubLya}, we obtain the
heteroclinic trajectory
represented in figure~\ref{fig:heteroclinic}. From now on, we will denote
this heteroclinic orbit by $\Het$. Note that this computation only
takes few seconds on a standard desktop computer.

\subsection{From One Orbit to Another}\label{sec:orbits}

Here, we construct two rather simple problems: first we compute an
optimal control using the Pontryagin Maximum Principle reaching the
heteroclinic orbit from the
Lyapunov orbit around $L_{1}$,  then we compute an optimal control to
reach the Lyapunov orbit around $L_{2}$ from the heteroclinic
orbit. This way, we get an admissible control that follows the null
control heteroclinic orbit during a
certain time.

\subsubsection{Around $L_1$}\label{sec:aroundL1}

\pardesc{Problem Statement.} Consider two points $\chi_{0}^{*}\in\Lya_{1}$
and $\chi_{1}^{*}\in\Het$, a time $t_{0}$ and an initial mass
$m_{0}^{*}=\SI{1500}{kg}$.
We apply the Pontryagin  Maximum Principle to the following
problem:\footnote{We assume that the reader is familiar with the
  principal concepts of the Pontryagin Maximum Principle. For
  details, see~\cite{Pontryagin, Trelat2008}.}

\begin{equation}
  \label{eq:OCPL1}
  \Prob{L_{1}}\left\{
    \begin{array}{l}
      \displaystyle \min \int_{0}^{t_{0}}\norm{u}{2}\D t,\\
      \displaystyle\dot x = F_{0}(x)+\frac{\epsilon}{m}\sum_{i=1}^{2}u_{i}F_{i}(x),\\
      \displaystyle\dot m = -\beta_{*}\epsilon\norm{u}{},\\
      \norm{u}{}\leq 1,\\
      x(0)=\chi_{0}^{*},\;m(0)=m_{0}^{*}\text{ and }x(t_{0})=\chi_{1}^{*}.
    \end{array}
  \right.
\end{equation}

Here, we have fixed the two points $\chi_{0}^{*}$ and
$\chi_{1}^{*}$ on the Lyapunov orbit and the heteroclinic orbit. We will
see how we choose these points later. We will
release the constraint on the position of these two points by an
optimization and satisfy the transversality
conditions for problem~\ref{eq:OCPgoal} in the last steps of our method.

Since the two points $\chi_{0}^{*}$ and $\chi_{1}^{*}$ belong to trajectories
with an energy greater than $\energy(L_{2})>\energy(L_{1})$, we know
that an admissible trajectory connecting $\chi_{0}^{*}$ to $\chi_{1}^{*}$ exists
(see~\cite{Daoud2011}).

If $t_{0}$ is greater than the \emph{minimum time}, we can show that
we are in the \emph{normal} case for the Pontryagin Maximum Principle,
that is to say $p^{0}$ can be normalized to $-1$
(see proposition~2 in~\cite{Caillau2012}). Although we have not proved
that this assumption holds, we will see that it is a reasonable one
because of the
construction of our two points. Moreover, because normality of
  the trajectories relies on the invariance of the target with respect
  to the zero control (see~\cite{Haberkorn2006}
  and~\cite{Chen}), the normality property holds for the
  \emph{targeted
    problem}~\eqref{eq:OCPgoal}.

We define the Hamiltonian as
$\Ham(x,m,p,p_{m},u)=-\norm{u}{2}+\prodscal{p}{F(x)}-\prodscal{p_{m}}{\beta_{*}\epsilon\norm{u}{}},$
where $F(x)=F_{0}(x)+\frac{\epsilon}{m}\sum_{i=1}^{2}u_{i}F_{i}(x)$,
$p\in\R^{4}$ and $p_{m}\in\R$.

To simplify the notation, we write:
\begin{equation*}
  \label{eq:Hamiltonian}
  \Ham(x,m,p,p_{m},u)=-\norm{u}{2}+H_{0}+H_{1}+H_{2}-\prodscal{p_{m}}{\beta_{*}\epsilon\norm{u}{}},
\end{equation*}
where $H_{i}=\prodscal{p}{F_{i}(x)}$, $i=0,\dots,2$.

Let us define $\varphi(p)=(p_{3},p_{4})$, thanks to the maximization
condition of the Pontryagin Maximum Principle, we get the optimal
control. Denoting by $y=(x,m,p,p_{m})$, let us
introduce the \emph{switching function}:
\begin{equation*}
  \label{eq:switch}
  \psi(y)=\frac{-\beta_{*}\epsilon p_{m}-\epsilon/m\norm{\varphi(p)}{}}{2}.
\end{equation*}
Then, the control is:
\begin{itemize}
\item if $\norm{\varphi(p)}{}\neq 0$, then
  \begin{equation*}
    \label{eq:control}
    \left\{
      \begin{array}{ll}
        u(y) = 0&\text{ if $\psi(y)\leq 0$},\\
        u(y) = \psi(y)\frac{\varphi(p)}{\norm{\varphi(p)}{}}&\text{ if $\psi(y)\in[0,1]$},\\
        u(y) = \frac{\varphi(p)}{\norm{\varphi(p)}{}}&\text{ else},\\
      \end{array}
    \right.
  \end{equation*}

\item if $\abs{\varphi(p)}=0$, then
  \begin{equation*}
    \label{eq:controlSphere}
    \left\{
      \begin{array}{ll}
        u(y) = 0&\text{ if $\psi(y)\leq 0$},\\
        u(y) \in \sphere(0,\psi(y))&\text{ if $\psi(y)\in[0,1]$},\\
        u(y) \in \sphere(0,1) &\text{ else},\\
      \end{array}
    \right.
  \end{equation*}
  where $\sphere(a,b)$ is the $\R^{2}$-sphere centered in $a$ with
  radius $b$. We will not take into account the singularity
    $\varphi(p)=0$. Hence, the control is
    continuous. This is one of the reasons why the numerical methods are
    easier for the minimization of the $L_{2}$-norm of $u$ than for
    the minimization of the $L_{1}$-norm.
\end{itemize}

In this problem, let us write the transversality conditions from the
Pontryagin Maximum Principle for the first problem~\eqref{eq:OCPL1}. The free
mass at the end of the transfer gives us : $p_{m}(t_{0})=0$. Moreover,
because of the final condition $x(t_{0})=\chi_{1}^{*}$, $p(t_{0})$ is
free. Finally, we are left to find $(p(0),p_{m}(0))$ such that the
final state
condition is satisfied.

We can write this problem as a shooting function. We denote
by $\flot^{\Ext}$ the extremal flow of the
extremal system. Hence, we define the shooting function:
\begin{equation}
  \label{eq:shootL1}
  \Ftir_{L_{1}}(p(0),p_{m}(0))
  = \begin{pmatrix}\flot^{\Ext}_{1,\dots,4}(\chi_{0}^{*},m_{0}^{*},p(0),p_{m}(0))-\chi_{1}\\[5pt]
    \flot^{\Ext}_{10}(\chi_{0}^{*},m_{0}^{*},p(0),p_{m}(0))
  \end{pmatrix}=
  \begin{pmatrix}
    \vect{0}\\[5pt] 0
  \end{pmatrix}.
\end{equation}

We compute the solution, that is to say $p(0)$ and $p_{m}(0)$ using
a shooting method (Newton-like method applied to~\eqref{eq:shootL1}). As
is well known, the
main difficulty is to \emph{initialize} the Newton-like algorithm. To do
this, we have used a continuation method. Let us explain the process.

\pardesc{Construction of $\chi_{0}^{*}$ and $\chi_{1}^{*}$.}
We want to realize the transfer from $\Lya_{1}$ to
$\Het$ and we have already computed the heteroclinic orbit. The method is
the following:
\begin{enumerate}
\item If we denote by $\chi_{\Het}^{L_{1}}$ the first point of the
  ``numerical'' heteroclinic orbit near the Lypunov orbit, we find $\chi_{\Lya_{1}}\in
  \Lya_{1}$ by minimizing the euclidean norm :
  $\chi_{\Lya_{1}}=\argmin_{\chi\in\Lya_{1}}\norm{\chi_{\Het}^{L_{1}}-\chi}{}$.
\item Then, we propagate backward in time $\chi_{\Lya_{1}}$ following the
  uncontrolled dynamics during a time $t_{\Lya_{1}}$ (smaller than the
  period of the Lyapunov orbit) to get $\chi_{0}^{*}$
\item We propagate forward in time $\chi_{\Het}^{L_{1}}$ during a reasonable
  time $t_{\Het}^{L_{1}}$ to get $\chi_{1}^{*}$ (small compared to the
  traveling time to reach the other extremity of $\Het$).
\end{enumerate}

We define the transfer time $t_{0}$ in~\eqref{eq:OCPL1} as
$t_{0}=t_{\Lya_{1}}+t_{\Het}^{L_{1}}.$

Although it seems to be a more simple problem than
problem~\eqref{eq:OCPgoal}, the main difficulty is
 still to initialize the shooting method. We
use a continuation method on the final state, using as a first simpler
problem a natural trajectory corresponding to a null control. Then step
by step, we reach the targeted  final point on the heteroclinic orbit,
as explained next.

\pardesc{Final State Continuation.} As explained in
section~\ref{sec:continuation}, we construct a family of problems
$\Prob{\lambda}$
depending continuously on one parameter $\lambda$ such that $\Prob{0}$
is easy to solve and $\Prob{1}$ corresponds to the targeted problem, that
is to say~\eqref{eq:OCPL1}.

First, let us define $\chi_{\Lya_{1}}^{\Nat}$ as the forward propagation of
$\chi_{\Lya_{1}}$  following the uncontrolled dynamics during
time~$t_{0}$. Then we define the family of problems:

\begin{equation*}
  \label{eq:OCPL1F}
  \Prob{L_{1}}^{\lambda}\left\{
    \begin{array}{l}
      \displaystyle \min \int_{0}^{t_{0}}\norm{u}{2}\D t,\\
      \displaystyle\dot x = F_{0}(x)+\frac{\epsilon}{m}\sum_{i=1}^{2}u_{i}F_{i}(x),\\
      \displaystyle\dot m = -\beta_{*}\epsilon\norm{u}{},\\
      \norm{u}{}\leq 1,\\
      x(0)=\chi_{0}^{*},\;m(0)=m_{0}^{*},\\[3pt]
      x(t_{0})=(1-\lambda)\chi_{\Lya_{1}}^{\Nat}+\lambda \chi_{1}^{*}.
    \end{array}
  \right.
\end{equation*}

\begin{figure}[ht]\centering
  \includegraphics[page=8]{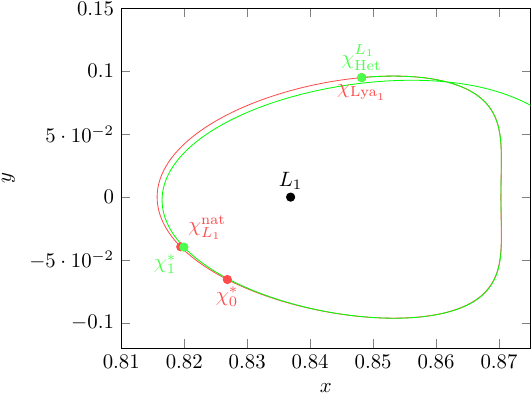}
  \caption{All relevant points in the construction of the problem.}
  \label{fig:PointsProbL1}
\end{figure}

Since $\chi_{\Lya_{1}}^{\Nat}$ corresponds to the uncontrolled dynamics,
the corresponding initial costate
$(p(0),p_{m}(0))$ is zero. Then, step by step, we initialize the
shooting method of $\Prob{L_{1}}^{\lambda_{i}}$ using the solution of
$\Prob{L_{1}}^{\lambda_{i-1}}$ to reach problem~\eqref{eq:OCPL1}. This
is done by a linear prediction, \latin{e.i.}, the solution of the two
previous iterations of the continuation are used to initialize the
resolution of the next step by a linear prediction.

Figure~\ref{fig:PointsProbL1} shows the different
points defined for some parameters described below.

\pardesc{Numerical Results.} We show here the numerical results for
this transfer. We choose a maximal thrust equal to \SI{60}{N}. We
postpone to  section~\ref{sec:multipleshooting} the problem of the
maximum thrust which should be very
small. In fact, a high thrust implies that the magnitude of the
costate stays very low, and it will be necessary for the multiple
shooting to converge. Indeed, for a non-saturating control, the
  higher the maximal magnitude of the thrust, the lower the control
  $u$ is between $[0,1]$, and so the lower the magnitude of $\psi(y)$
  and thus of the costate.
Moreover, we choose the two times of propagation in the normalized
system as
$t_{\Lya_{1}}=1.0$, and $t_{\Het}^{L_{1}}=2.0.$

\begin{figure}[ht]\centering
  \includegraphics[page=9]{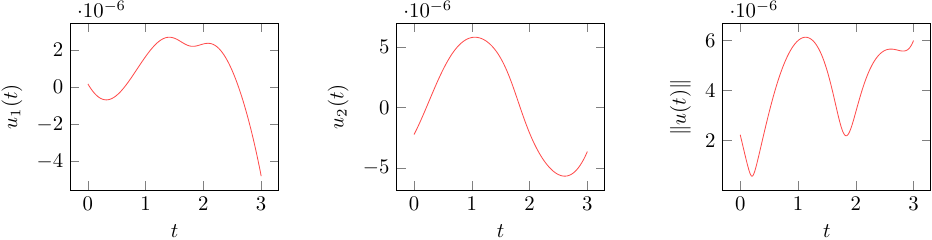}
  \caption{Command to realize the optimal transfer from the Lyapunov
    orbit to the heteroclinic orbit. We plot $u(\cdot)\leq 1$ as a
    function of the normalized time.}
  \label{fig:commandeL1}
\end{figure}

We obtain the optimal trajectory plotted in
figure~\ref{fig:PointsProbL1}. The optimal command is
shown in figure~\ref{fig:commandeL1}. One can see that we are
far from the saturation of the command, indeed, the maximum value is
approximately \si{6e-06}, whereas we are constrained by one. We postpone the
discussion on the real value in Newton to the final trajectory.

This continuation gives us an initial adjoint vector (costate) that we
will denote
 by $p_{0}^{*}$ and $p_{m}^{0*}$ in the remainder of this work.

 \subsubsection{Around $L_2$}\label{sec:aroundL2}

We design a very similar problem around $L_{2}$.

\pardesc{Problem Statement.} Consider two points $\chi_{2}^{*}\in\Het$
and $\chi_{3}^{*}\in\Lya_{2}$, a time $t_{2}$ and an initial mass
$m_{2}^{*}$. \footnote{We will understand why we use 2 as subscript.}
The mass $m_{2}^{*}$ is the final mass obtained after  solving for
the transfer around $L_{1}$ (between the two problem we follow a
heteroclinic orbit without any fuel consumption).
We apply  the Pontryagin  Maximum Principle to the following problem:

\begin{equation}
  \label{eq:OCPL2}
  \Prob{L_{2}}\left\{
    \begin{array}{l}
      \displaystyle \min \int_{0}^{t_{2}}\norm{u}{2}\D t,\\
      \displaystyle\dot x = F_{0}(x)+\frac{\epsilon}{m}\sum_{i=1}^{2}u_{i}F_{i}(x),\\
      \displaystyle\dot m = -\beta_{*}\epsilon\norm{u}{},\\
      \norm{u}{}\leq 1,\\
      x(0)=\chi_{2}^{*},\;m(0)=m_{2}^{*}\text{ and }x(t_{2})=\chi_{3}^{*}.
    \end{array}
  \right.
\end{equation}

As before, we have fixed $\chi_{2}^{*}$ and $\chi_{3}^{*}$ on the
heteroclinic and Lyapunov orbits. The final steps will allow us to
release these constraints.

Since the problem is very similar to the problem around $L_{1}$, we
have the same Hamiltonian and the same expression of the control
$u$. Hence, we get the following shooting function:

\begin{equation*}
  \label{eq:shootL2}
  \Ftir_{L_{2}}(p(0),p_{m}(0))
  = \begin{pmatrix}\flot^{\Ext}_{1,\dots,4}(\chi_{2}^{*},m_{2}^{*},p(0),p_{m}(0))-\chi_{3}^{*}\\[5pt]
    \flot^{\Ext}_{10}(\chi_{2}^{*},m_{2}^{*},p(0),p_{m}(0))
  \end{pmatrix}=
  \begin{pmatrix}
    \vect{0}\\[5pt] 0
  \end{pmatrix}.
\end{equation*}

\pardesc{Construction of $\chi_{2}^{*}$ and $\chi_{3}^{*}$.} We construct the two
points following the same method.

\begin{enumerate}
\item If we denote by $\chi_{\Het}^{L_{2}}$ the last point of the
  heteroclinic orbit near the Lyapunov orbit, we find $\chi_{\Lya_{2}}\in
  \Lya_{2}$ minimizing the euclidean norm :
  $\chi_{\Lya_{2}}=\argmin_{\chi\in\Lya_{2}}\norm{\chi_{\Het}^{L_{2}}-\chi}{}$.
\item Then, we propagate forward $\chi_{\Lya_{2}}$ following the
  uncontrolled dynamics during a time $t_{\Lya_{2}}$ (smaller than the
  period of Lyapunov orbit) to get $\chi_{3}^{*}$.
\item We propagate backward the $\chi_{\Het}^{L_{2}}$ during a reasonable
  time $t_{\Het}^{L_{2}}$ to get $\chi_{2}^{*}$ (small compared to the
  traveling time to reach the other extremity).
\end{enumerate}

We define the transfer time $t_{2}$ in~\eqref{eq:OCPL2} as
$t_{2}=t_{\Lya_{2}}+t_{\Het}^{L_{2}}.$

\pardesc{Final State Continuation.} As before, we construct a family
of problems
$\Prob{\lambda}$
depending continuously on one parameter $\lambda$ such that $\Prob{0}$
is easy to solve and $\Prob{1}$ corresponds to the targeted problem, that
is to say~\eqref{eq:OCPL1}.

First, let us define $\chi_{\Het}^{\Nat}$ as the forward propagation of
$\chi_{\Het}$  following the uncontrolled dynamics during the
time~$t_{2}$. Then we define the family of problems:

\begin{equation*}
  \label{eq:OCPL2F}
  \Prob{L_{2}}^{\lambda}\left\{
    \begin{array}{l}
      \displaystyle \min \int_{0}^{t_{2}}\norm{u}{2}\D t,\\
      \displaystyle\dot x = F_{0}(x)+\frac{\epsilon}{m}\sum_{i=1}^{2}u_{i}F_{i}(x),\\
      \displaystyle\dot m = -\beta_{*}\epsilon\norm{u}{},\\
      \norm{u}{}\leq 1,\\
      x(0)=\chi_{2}^{*},\;m(0)=m_{2}^{*},\\[3pt]
      x(t_{2})=(1-\lambda)\chi_{\Het}^{\Nat}+\lambda \chi_{3}^{*}.
    \end{array}
  \right.
\end{equation*}

\pardesc{Numerical Results.} As before, we set $\Tmax=\SI{60}{N}$,
and we compute the continuation for the two times chosen as
$t_{\Het}^{L_{2}}=2.0,$ and $t_{\Lya_{2}}=1.0$.

\begin{figure}[ht]\centering
  \includegraphics[page=10]{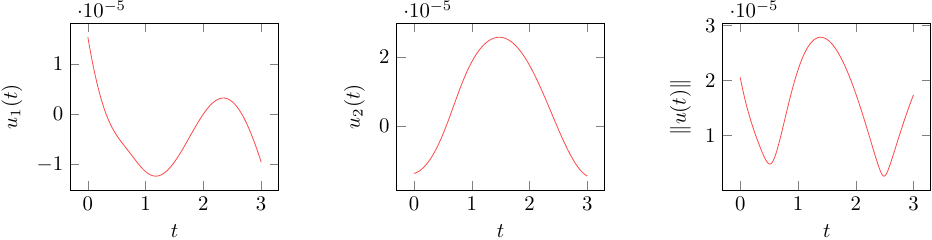}
  \caption{Command to realize the optimal transfer from the Lyapunov
    orbit to the heteroclinic orbit. We plot $u(\cdot)\leq 1$
    as a function of the normalized time.}
  \label{fig:commandeL2}
\end{figure}

Figure~\ref{fig:commandeL2} shows the optimal control to realize
the final transfer
from the heteroclinic orbit to the Lyapunov one around $L_{2}$. We see
that, once again, we are far from the saturation of $u$.

This continuation gives us an initial costate that we will denote by
$p_{2}^{*}$ and $p_{m}^{2*}$ in the remainder of this work.
\bigskip

Table~\ref{tab:resultsLi} sums up all the parameters for the
continuation computation. We observe that, because we are using indirect
shooting methods, the computation is very fast even though it is performed on
a simple desktop computer or on a single-board computer (the
Raspberry Pi).

\begin{table}[ht]
  \centering
  \begin{tabular}{cccccc}
    \toprule
     $\Isp$ & $\G$ & Earth Mass & Moon Mass & Distance & Period  \\
     \midrule
    \SI{2000}{s} & \SI{9.8}{m^3 kg^{-1} s^{-2}} & \SI{5.972E24}{kg} &
    \SI{7.349E22}{kg} &  \SI{384402e3}{m}  & \SI{2.361e6}{s}\\
    \bottomrule
  \end{tabular}
  \bigskip

  \begin{tabular}{cccc}
    \toprule
    Transfer & Iterations &  Cost & $\Tmax$ \\
    \midrule
   $L_{1}$ & 21 &  \num{6.30967e-11}& \SI{60}{N}\\
   $L_{2}$ & 19  &\num{9.06124e-10}& \SI{60}{N}\\
   \bottomrule
 \end{tabular}
 \bigskip

 \begin{tabular}{ccc}
   \toprule
   System & Transfer &  Execution time\\
   \midrule
   Core i7 & $L_{1}$ & \texttt{98\% cpu 2,821s
     total} \\
   & $L_{2}$ &\texttt{96\% cpu 1,439s total}\\
   Raspberry Pi A & $L_{1}$ & \texttt{38\% cpu
     8,009s total}\\
   & $L_{2}$ & \texttt{22\% cpu 7,879s total}\\
   \bottomrule
 \end{tabular}

  \caption{Numerical results for the two transfers around $L_1$ and
    $L_2$. Computations are performed on a simple laptop Core i7, and on a
    Raspberry Pi A, a credit card-sized
    single-board computer.}
  \label{tab:resultsLi}
\end{table}

\subsection{Multiple Shooting}\label{sec:multipleshooting}

Thanks to the results from previous sections, we have  designed  an
admissible
control to
perform  the transfer from a Lyapunov orbit around $L_{1}$ to a
Lyapunov orbit around $L_{2}$. We first reach a point on a
heteroclinic orbit, then we follow the natural dynamics (null
control), and finally reach a point on the final Lyapunov orbit from a
certain point on the heteroclinic orbit. This admissible trajectory is
however not energy optimal,
since the stay on the heteroclinic orbit is forced.

These two points on the heteroclinic orbit were arbitrarily
chosen. There is no guarantee that they provide a good choice in
terms of
optimality. Hence, we want to release the constraints on the
position of these two points. We
use a multiple shooting method on top of the first two local transfer
to get a better optimum.

\begin{figure}[ht]\centering
  \includegraphics[page=11]{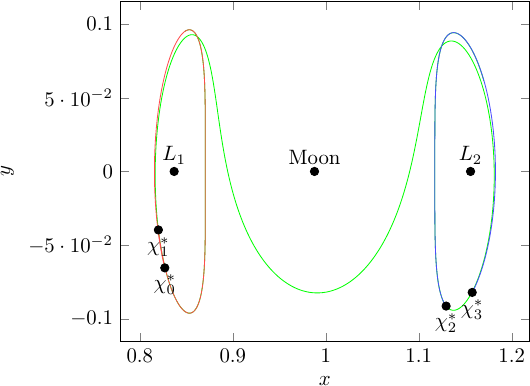}
  \caption{Admissible trajectory in three parts.}
  \label{fig:multiple}
\end{figure}

Let us describe how we state the multiple shooting problem.
As we can see in figure~\ref{fig:multiple}, there are two points
$\chi_{1}^{*}$ and $\chi_{2}^{*}$ belonging to the heteroclinic orbit that we want
to free. Moreover, we have three times:
\begin{itemize}
\item $t_{0}$ which is the time defined for the transfer around $L_{1}$;
\item $t_{2}$  which is the time defined for the transfer around
  $L_{2}$;
\item $t_{1}$ which is the total time of the computed heteroclinic
  orbit minus the two times $t_{\Het}^{L_{1}}$ and $t_{\Het}^{L_{2}}$
  used in the two previous transfers.
\end{itemize}
We define $t_{\tot}=t_{0}+t_{1}+t_{2}$ and we
write a new optimal control problem with the same structure as the
previous one around $L_{1}$ and $L_{2}$.

\begin{equation}
  \label{eq:OCPTOT}
  \Prob{\tot}\left\{
    \begin{array}{l}
      \displaystyle \cost_{\tot}=\min \int_{0}^{t_{\tot}}\norm{u}{2}\D t,\\
      \displaystyle\dot x = F_{0}(x)+\frac{\epsilon}{m}\sum_{i=1}^{2}u_{i}F_{i}(x),\\
      \displaystyle\dot m = -\beta_{*}\epsilon\norm{u}{},\\
      \norm{u}{}\leq 1,\\
      x(0)=\chi_{0}^{*}\in\Lya_{1},\;m(0)=m_{0}^{*},\\[3pt]
      x(t_{\tot})=\chi_{3}^{*}\in\Lya_{2}.
    \end{array}
  \right.
\end{equation}

As before, we apply the Pontryagin  Maximum Principle to get a
necessary condition for the optimal control. We are able
to write the control $u$ with respect to the state $(x,m)$ and the costate
$(p,p_{m})$, we can write a shooting function,  with the same results
as the ones obtained in section~\ref{sec:orbits}.

Thanks to the following method, we get an admissible trajectory in
three parts, and it is quite natural to use it to construct a multiple
shooting function. We define
\[Z=(\underbrace{p_{0},p_{m}^{0}}_{P_{0}},\underbrace{\chi_{1},m_{1}}_{X_{1}},
\underbrace{p_{1},p_{m}^{1}}_{P_{1}},\underbrace{\chi_{2},m_{2}}_{X_{2}},\underbrace{p_{2},p_{m}^{2}}_{P_{2}})\in\R^{25},\]
then we write the multiple shooting function with two matching
conditions on the state and the costate, the final state condition, and the
free final mass:
\begin{equation}
  \label{eq:multipleshooting}
  \Ftir_{\multi}(Z)=
  \begin{pmatrix}
    \flot^{\Ext}_{1,\dots,5}(\chi_{0}^{*},m_{0}^{*},P_{0})-X_{1}\\
    \flot^{\Ext}_{6,\dots,10}(\chi_{0}^{*},m_{0}^{*},P_{0})-P_{1}\\
    \flot^{\Ext}_{1,\dots,5}(X_{1},P_{1})-X_{2}\\
    \flot^{\Ext}_{6,\dots,10}(X_{1},P_{1})-P_{2}\\
    \flot^{\Ext}_{1,\dots,4}(X_{2},P_{2})-\chi_{3}^{*}\\
    \flot^{\Ext}_{10}(X_{2},P_{2})\\
  \end{pmatrix}.
\end{equation}

We want to find the vector $Z$ such that
$\Ftir_{\multi}(Z)=0$, and as in previous sections, we use a Newton-like
algorithm. The main difficulty is as usual to initialize the
algorithm. This time,  it is done by the previous
\emph{local} transfers, since we chose:\footnote{Note that the
  notation $p_{1}$ and $p_{2}$ is not for the first and second
  components of the costate but for two different costate belonging
  to $\R^{4}$.}
\begin{equation*}
  \label{eq:initMultiple}
  \left\{
  \begin{array}{lllll}
    p_{0}=p_{0}^{*},&  \chi_{1}=\chi_{1}^{*}, & p_{1}= \vect{0},& \chi_{2}=\chi_{2}^{*},&p_{2}=p_{2}^{*},\\
    p_{m}^{0}=p_{m}^{0*},&
    m_{1}=m_{2}^{*},&
    p_{m}^{1}=0,&
    m_{2}=m_{2}^{*},&
    p_{m}^{2}= p_{m}^{2*}.
  \end{array}
  \right.
\end{equation*}
The choices  $m_{1}=m_{2}^{*}$, $p_{1}=\vect{0}$ and $p_{m}^{1}=0$ are
made
because we initialize the trajectory with a heteroclinic part, that is
to say with a null control and without consumption of mass.

The Newton-like algorithm gives us a complete trajectory which is not
constrained to follow the heteroclinic orbit. In
figure~\ref{fig:multipleTraj} and figure~\ref{fig:multipleCom}, we can
see the trajectory and the associated control.

We keep the maximum thrust equal to \SI{60}{N} to allow the
Newton-like algorithm to converge. But, we want to be able to give
the right specification for the engine of the spacecraft. Let us see how
we make this possible.

\subsection{Thrust Continuation}

Using, the continuation method we want to constrain the thrust to a
real value for a low-thrust engine, let us say \SI{0.3}{N}. To do
that, we construct a family of problems as before. Let us denote by
$\epsilon_{0}$ the initial maximal thrust in normalized units corresponding to
$\Tmax=\SI{60}{N}$. Similarly, let us denote by $\epsilon_{1}$ the
maximal thrust that we want to get corresponding to
$\Tmax=\SI{0.3}{N}$. Finally, we define the
maximal continuation thrust:
\[\epsilon_{\lambda}=(1-\lambda)\epsilon_{0}+\lambda\epsilon_{1}.\]
We can now define the family of problems:

\begin{equation*}
  \label{eq:OCPThrust}
  \Prob{\thrust}^{\lambda}\left\{
    \begin{array}{l}
      \displaystyle \min \int_{0}^{t_{\tot}}\norm{u}{2}\D t,\\
      \displaystyle\dot x = F_{0}(x)+\frac{\epsilon_{\lambda}}{m}\sum_{i=1}^{2}u_{i}F_{i}(x),\\
      \displaystyle\dot m = -\beta_{*}\epsilon_{\lambda}\norm{u}{},\\
      \norm{u}{}\leq 1,\\
      x(0)=\chi_{0}^{*}\in\Lya_{1},\;m(0)=m_{0}^{*},\\
      x(t_{\tot})=\chi_{3}^{*}\in\Lya_{2}.
    \end{array}
  \right.
\end{equation*}

We solve each step of the continuation with the previously defined
multiple shooting method~\eqref{eq:multipleshooting}. This way, we
manage to constrain the
thrust to the given engine value. Since the control is smaller than
\SI{0.3}{N}, this continuation is easy, and the command does
not
change during it. In section~\ref{sec:numRes} we summarize the
numerical results.
Let us remark that in our numerical experiment, the continuation is
done with 22 iterations.

\subsection{Optimization of the Terminal Points}

The last remaining step is to free the initial and final
points. The only constraints  are that $x(0)$
has to belong to $\Lya_{1}$ and $x(t_{\tot})$ to $\Lya_{2}$.
To simplify the problem, we have fixed by construction two points
$\chi_{0}^{*}$ on $\Lya_{1}$ and $\chi_{3}^{*}$ on $\Lya_{2}$. Now we
want to find the optimal points on these two periodic orbits.
So we want to
solve the very general problem~\eqref{eq:OCPgoal}. The
Pontryagin Maximum Principle gives us two transversality conditions
that we have to satisfy:
\begin{equation}\label{eq:transverse}
  p_{1,\dots,4}(0)\ \bot\ T_{x(0)}\Lya_{1}\quad
  \text{and}\quad p_{1,\dots,4}(t_{\tot})\ \bot\ T_{x(t_{\tot})}\Lya_{2},
\end{equation}
where the notation $T_{x}M$ stands for the usual tangent space to
$M$ at the point $x$ (these conditions can be written as soon as
the tangent space is well defined).

To perform this optimization we consider the two previously chosen points
$\chi_{0}^{*}\in\Lya_{1}$ and
$\chi_{3}^{*}\in\Lya_{2}$. First we perturb the point around
$\Lya_{2}$ following the decrease of the transversality
condition until it changes sign so to find a good zero for
the transversality
condition. Since we checked that the evolution of this
transversality
condition along the periodic orbit is not monotone, we are just able
to reach a local minimum. By doing this we manage to reach a
transversality condition at $t_{\tot}$ around $\num{1E-08}$. Secondly, we
realize
the same perturbation along $\Lya_{1}$ and we manage to reach a value
around $\num{1E-08}$. We have checked that the inverse process beginning
with the point on $\Lya_{1}$ gives the same result.

Although this seems to cause very little change on the transfer (see
numerical results in the next section), the
structure of the control is completely changed. We will describe this
result in depth in the next
section.

\begin{Remark}
To perform this optimization, we could use a gradient method on the
one-dimensional periodic orbits initializing it with the solution of
 problem~\eqref{eq:OCPTOT}.
\end{Remark}

\subsection{Numerical Results}\label{sec:numRes}

\begin{figure}[ht]\centering
  \includegraphics[page=12]{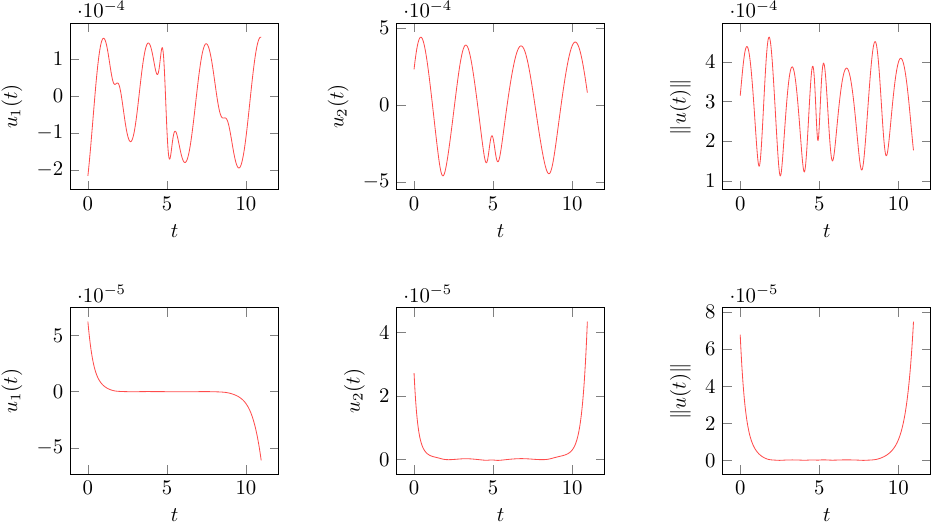}
  \caption{Command to realize the optimal transfer from the Lyapunov
    orbit around $L_{1}$ to the Lyapunov orbit around $L_{2}$. We plot
    $u(\cdot)\leq 1$ before
    the last optimization step on the first row (we chose two points
    on $\Lya_{1}$ and $\Lya_{2}$) and after the last optimization step
    consisting in getting the general transversality conditions
    (second row). We
    can observe the good turnpike property of the second control.}
  \label{fig:multipleCom}
\end{figure}

Recall that we use the CRTBP parameters given in
table~\ref{tab:resultsLi}. We observe in figure~\ref{fig:multipleCom} that
the last
optimization step changes the shape of the control. Indeed, by
construction, we make the spacecraft go onto the heteroclinic
orbit before we free that constraint. Hence, it can be expected that
the mission has \emph{turnpike
  properties} (see~\cite{trelatzuazua}). That is to say the optimal
solution settled in large time consists
approximately of three pieces, the first
and the last of which being transient short-time arcs, and the middle
piece being a long-time arc staying
exponentially close to the optimal steady-state solution. In
figure~\ref{fig:multipleCom}, we see that before  the
transversality conditions are satisfied following the last
optimization step, the command
structure does not have the shape of a turnpike command. Control is
spread along the trajectory. After the last optimization step, the
control is clearly a turnpike control and the trajectory consists
approximately in three pieces as expected.

\begin{figure}[ht]\centering
  \includegraphics[page=13]{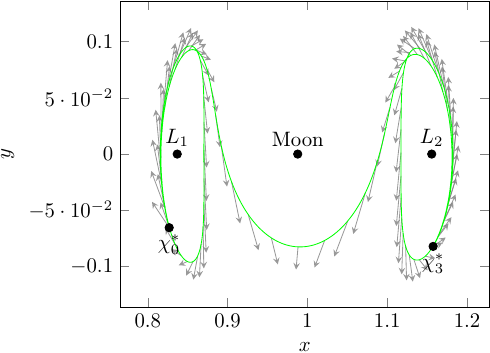}
  \hfill
  \includegraphics[page=14]{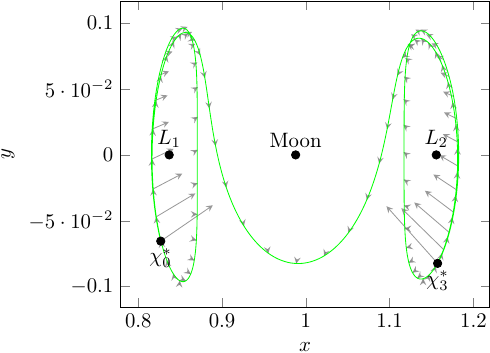}
  \caption{Optimal trajectory. On the left the optimal trajectory
    with $\chi_{0}^{*}$ and $\chi_{3}^{*}$ \emph{fixed} on $\Lya_{1}$ and
    $\Lya_{2}$. On the right, the optimal trajectory with
    $\chi_{0}^{*}$ and $\chi_{3}^{*}$ \emph{free} on $\Lya_{1}$ and
    $\Lya_{2}$. The control is represented by arrows.}
  \label{fig:multipleTraj}
\end{figure}

We show  in figure~\ref{fig:multipleTraj} the two corresponding
trajectories. We observe that, to satisfy to transversality conditions
corresponding to $x(0)\in\Lya_{1}$ and
$x(t_{\tot})\in\Lya_{2}$,  the two fixed points were note moved very
much.

\pardesc{Cost.} In this problem we are minimizing the cost
$\int_{0}^{t_{\tot}}\norm{u}{2}\D t$.
We consider a mass evolving dynamical system, and a maximum
thrust so to try to compare fairly the cost with other results, we
define three different costs:
\begin{equation}
  \label{eq:costs}
  \cost_{\tot}^{1}=\int_{0}^{t_{\tot}}\norm{u(t)}{2}\D t,\quad
  \cost_{\tot}^{2}=\int_{0}^{t_{\tot}}\dfrac{\epsilon^{2}}{m^{2}(t)}\norm{u(t)}{2}\D
  t,\;\text{ and }\;\cost_{\tot}^{3}=\int_{0}^{t_{\tot}}\dfrac{\Tmax^{2}}{m^{2}(t)}\norm{u(t)}{2}\D
  t.
\end{equation}

Results are summarized in table~\ref{tab:final}. We observe that
whereas the two points $\chi_{0}^{*}$ and $\chi_{3}^{*}$ are not
perturbed very much to satisfy the general transversality conditions, for
the costs and the mass consumption, it is really an improvement.

\begin{table}[ht]
  \centering

  \begin{tabular}{ccc}
    \toprule
    Initial Mass &  Transfer time & $\Tmax$\\
    \midrule
    \SI{1500}{kg} & \num{10.96139} or
    \SI{47.67}{days} & \SI{0.3}{N}\\
    \bottomrule
  \end{tabular}
  \bigskip

  \begin{tabular}{ccccc}
    \toprule
    &$\cost_{\tot}^{1}$& $\cost_{\tot}^{2}$& $\cost_{\tot}^{3}$  & Mass of fuel\\
    \midrule
 Problem~\eqref{eq:OCPTOT}&   \num{1.0650187e-06} & \num{ 5.7479872e-09}& \num{1.8527847e-13}& \SI{0.0186878}{kg} \\[5pt]
 Problem~\eqref{eq:OCPgoal}&   \num{2.2305967e-09} & \num{ 1.2038555e-11}& \num{ 3.8804630e-16}& \SI{3.6709589e-04}{kg} \\
 \bottomrule
  \end{tabular}
  \bigskip

  \begin{tabular}{ccc}
    \toprule
    &System  &  Execution time\\
    \midrule
    Problem~\eqref{eq:OCPTOT}  & Core i7   & \texttt{99\% cpu 26,912s total} \\
    Problem~\eqref{eq:OCPgoal}& Core i7 &\texttt{99\% cpu 1min18,64s total }\\
    Problem~\eqref{eq:OCPTOT} &Raspberry Pi A & \texttt{20\% cpu
      15min3,545s total}\\
    Problem~\eqref{eq:OCPgoal}& Raspberry Pi A & \texttt{23\% cpu
      56min45,921s total}\\
    \bottomrule
  \end{tabular}

  \caption{Numerical results for the final trajectory of the first
    mission obtained after
    the multiple shooting with fixed departure and final points
    (problem~\eqref{eq:OCPTOT}) and for the optimized departure and final
    points on $\Lya_{1}$ and $\Lya_{2}$ (problem~\eqref{eq:OCPgoal}).}
  \label{tab:final}
\end{table}

\section{Variant of the mission}

In this section, we show two other applications of our method to
design two different missions. The first one is the Lyapunov to
Lyapunov mission but with different energies and a heteroclinic orbit
with two revolutions around the Moon.

The second mission is a Halo to Halo mission with two different
energies and with no heteroclinic orbit. In this case, we use two
trajectories belonging to two invariant manifolds.

\subsection{Heteroclinic Orbit with Two Revolutions}\label{sec:mission2}

In this section we present another mission going from a Lyapunov orbit
around $L_{1}$ to a Lyapunov orbit around $L_{2}$. We follow exactly
the same method at the one we presented except that we find
the second intersection of manifolds (instead of the first) and  we compute
the second crossing through the plane $U_{2}$ on both sides with the
stable and unstable manifolds. Our final trajectory will perform two
revolutions
around the Moon. Because the $\Lya_{2}$ target is invariant
  with respect to the zero control, the larger the duration of the
  heteroclinic orbit, the smaller (and better) the fuel consumption. Considering
  that, we expect  a better cost for the transfer.

\subsubsection{The Heteroclinic orbit}

\begin{figure}[ht]\centering
  \includegraphics[page=15]{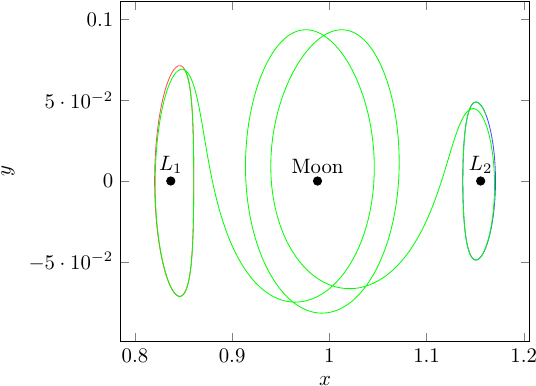}
  \caption{Heteroclinic orbit between two Lyapunov orbits in the
    Earth-Moon system. We get a travel time of $11.699681461946$
    (normalized time) or $50.883$ days.}
  \label{fig:hetero2}
\end{figure}

To compute the heteroclinic orbit with two revolutions around the
Moon, we
have to choose a certain energy allowing the second intersection to
exist. We have chosen the energy (we follow \cite{Epenoy} to motivate
this choice) $\energy_{\Lya_{1,2}} = -1.5890,$
and computed the heteroclinic orbit plotted in
figure~\ref{fig:hetero2}. There are indeed two revolutions around the
Moon.

\subsubsection{Two Local Transfers}

As before, we compute two local transfers. One from the periodic orbit
around $L_{1}$ to the heteroclinic orbit, and another from the end of
the heteroclinic orbit to the periodic orbit around $L_{2}$. We choose
the maximal thrust equal to $\SI{60}{N}$ as before to help the success
of the shooting. We do not
report the partial results here as they are comparable to the ones of
the previous mission. Thanks to this step, we obtain an
admissible trajectory in three parts, one controlled to reach the
heteroclinic orbit (the \emph{turnpike}), the second part is
the uncontrolled heteroclinic orbit, and the last part is a
controlled one from the heteroclinic to the Lyapunov orbit around
$L_{2}$.

\subsubsection{Multiple Shooting Method}

As before, to free the two matching connections on the heteroclinic
orbit and to decrease the maximum thrust we use a multiple shooting
method associated with a continuation method. Since the transfer
time is larger than for the previous
mission, we have to add some grid points along the heteroclinic orbit
(which are initialized with a null adjoint vector). This is due to the
very unstable nature of the hamiltonian system.  Here we chose 5 grid
points. Thanks to the multiple shooting method and a thrust
continuation, we manage to reach the
required maximal thrust: $\Tmax=\SI{0.3}{N}$ and we get an admissible
trajectory with two fixed points on $\Lya_{1}$ and $\Lya_{2}$. The
last step consists in finding the optimal departure and arrival points
on the two periodic orbits.

\subsubsection{Optimization of the Terminal Points}

Once again, because we have simplified the problem by fixing the
departure and arrival points on $\Lya_{1}$ and $\Lya_{2}$, we want to
free these points on the periodic orbits to satisfy the general
transversality conditions~\eqref{eq:transverse}. As before, we
perturb first $\chi_{3}^{*}\in\Lya_{2}$ following the decrease of the
transversality condition and we do the same with $\chi_{0}^{*}$. We
manage to satisfy the transversality conditions up to $\num{1E-09}$.

\subsubsection{Results}

\begin{figure}[ht]\centering
  \includegraphics[page=16]{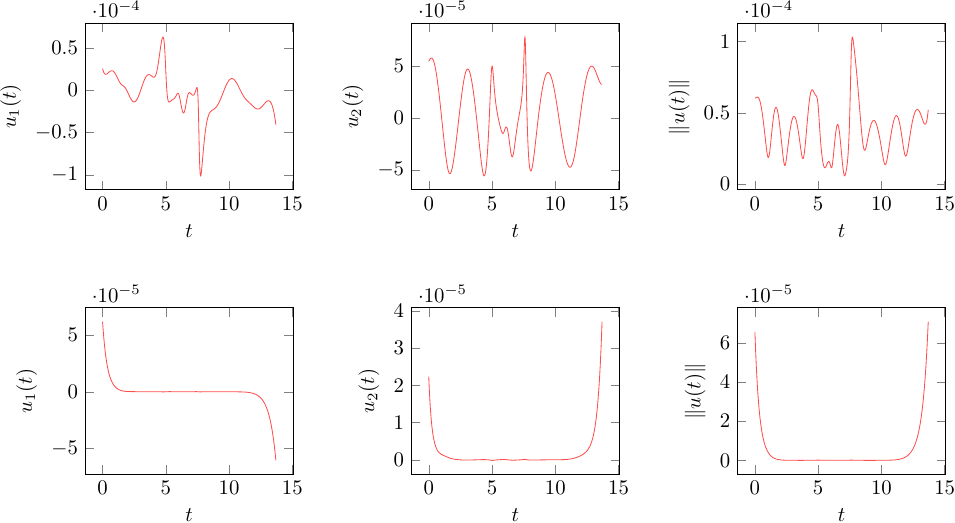}
  \caption{Command to realize the optimal transfer from the Lyapunov
    orbit to the heteroclinic orbit. We plot $u(\cdot)\leq 1$ before
    the last optimization step on the first row (we chose two points
    on $\Lya_{1}$ and $\Lya_{2}$) and after the last optimization step
    consisting in getting the general transversality conditions. We
    can observe the good turnpike property of the second control.}
  \label{fig:command2}
\end{figure}

We plot in figure~\ref{fig:command2} the command before and after the
last optimization step. We
observe the same phenomenon as for the previous mission. Indeed, before
we satisfy the transversality conditions, the command does not have the
\emph{turnpike structure}, that is to say, the three parts, first a
short thrust to reach the \emph{highway} (or turnpike), then a null controlled part,
and finally a controlled part to reach the periodic orbit.

\begin{figure}[ht]\centering
  \includegraphics[page=17]{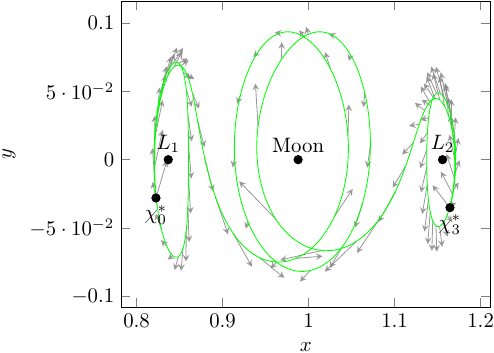}
  \hfill
  \includegraphics[page=18]{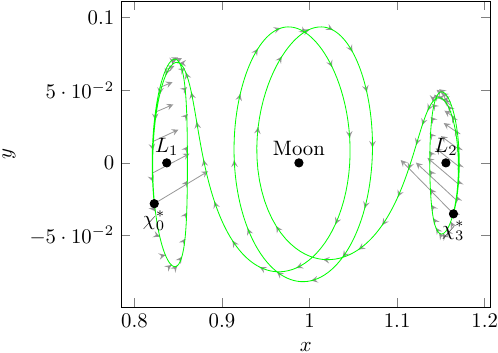}
  \caption{Optimal trajectory. On the left the optimal trajectory
    with $\chi_{0}^{*}$ and $\chi_{3}^{*}$ \emph{fixed} on $\Lya_{1}$ and
    $\Lya_{2}$. On the right, the optimal trajectory with
    $\chi_{0}^{*}$ and $\chi_{3}^{*}$ \emph{free} on $\Lya_{1}$ and
    $\Lya_{2}$. The control is represented by arrows.}
  \label{fig:traj2}
\end{figure}

Whereas the perturbations of the two points $\chi_{0}^{*}$ and
$\chi_{3}^{*}$ to satisfy the transversality conditions are very small
(see figure~\ref{fig:traj2}), the structure of the control is very
different and the costs are much smaller after getting the
transversality conditions. We summarize the numerical results in
table~\ref{tab:results2}.

\begin{table}[ht]
  \centering

  \begin{tabular}{ccc}
    \toprule
    Initial Mass &  Transfer time & $\Tmax$\\
    \midrule
    \SI{1500}{kg} & \num{13.699681461} or
    \SI{59.582}{days} & \SI{0.3}{N}\\
    \bottomrule
  \end{tabular}
  \bigskip

  \begin{tabular}{ccccc}
    \toprule
    &$\cost_{\tot}^{1}$& $\cost_{\tot}^{2}$& $\cost_{\tot}^{3}$  & Mass of fuel\\
    \midrule
    Problem~\eqref{eq:OCPTOT}&   \num{ 2.4638905e-08} & \num{1.3297667e-10}& \num{4.2863204e-15}& \SI{0.0030131}{kg} \\[5pt]
    Problem~\eqref{eq:OCPgoal}&   \num{ 1.9695934e-09} & \num{ 1.0629917e-11}& \num{ 3.4264079e-16}& \SI{3.3599750e-04}{kg} \\
    \bottomrule
  \end{tabular}
  \bigskip

  \begin{tabular}{ccc}
    \toprule
    &System  &  Execution time\\
    \midrule
    Problem~\eqref{eq:OCPTOT}  & Core i7   & \texttt{99\% cpu 44,949s total} \\
    Problem~\eqref{eq:OCPgoal}& Core i7 &\texttt{99\% cpu 2min54,79s total}\\
    Problem~\eqref{eq:OCPTOT} &Raspberry Pi A & \texttt{33\% cpu
      22min52,8s total}\\
    Problem~\eqref{eq:OCPgoal}& Raspberry Pi A & \texttt{29\% cpu
      1h32min46s total}\\
    \bottomrule
  \end{tabular}

  \caption{Numerical results for the final trajectory of the second
    mission obtained after
    the multiple shooting with fixed departure and final points
    (problem~\eqref{eq:OCPTOT}) and for the optimized departure and final
    points on $\Lya_{1}$ and $\Lya_{2}$ (problem~\eqref{eq:OCPgoal}).}
  \label{tab:results2}
\end{table}

\subsection{Halo to Halo Mission}\label{sec:HaloHalo}

In this section, we will adapt the previous method  to
another mission: a Halo to Halo mission. Halo orbits are periodic
orbits around equilibrium points like Lyapunov orbits but in the
spatial dynamics. Because of that, we consider in this section all the
previous concepts and results presented in the
sections~\ref{sec:mission} and~\ref{sec:dynamics} extended to the
spatial
configuration.

For the Halo to Halo mission, because we are in the spatial case and for the
energies of the periodic orbits that we have chosen, the intersection
between unstable and stable manifolds does not exist. However, our method is
still valid and can be applied.

We will first design an admissible
trajectory with 5 parts:
\begin{enumerate}
  \item first, we propagate the unstable and stable manifolds from
    $L_{1}$ and $L_{2}$ as described in section~\ref{sec:het}. We
    compute, in the plane $U_{2}$, the two points (one on each
    manifolds) that minimize the distance in position and velocity. This
    gives us two trajectories.
\item Then, we compute the optimal transfer from a fixed point on the
  Halo orbit around $L_{1}$ to
  a fixed point on the trajectory on the associated unstable manifold.
\item We compute a transfer from a fixed point on the trajectory
  of the unstable manifold from the Halo orbit around $L_{1}$ to a
  fixed point on the
  trajectory on the stable manifold of the Halo orbit around $L_{2}$.
\item We then compute the optimal control to reach a fixed point on
  the Halo orbit around $L_{2}$ from a fixed point on the trajectory
  of the associated stable manifold.
\end{enumerate}

 With this admissible trajectory in 5 parts (with two
  uncontrolled parts), we initialize a multiple shooting method to get
  an optimal trajectory reaching a fixed point on the Halo orbit around
  $L_{2}$ from a fixed point on the Halo orbit around
  $L_{1}$. Finally, following the method described for the Lyapunov to
  Lyapunov mission, we optimize the position of the end points.

As we can see from this example, the method is quite general, and we
can think about applying it for much more complex missions designed by
patching together ``manifold'' parts.

\subsubsection{Free Parts on Manifolds}

As described in the introduction, we will compute two trajectories on
unstable and stable manifolds respectively from the Halo orbit around
$L_{1}$ and from the Halo orbit around $L_{2}$.

\begin{figure}
  \centering
  \includegraphics[page=19, scale=0.9]{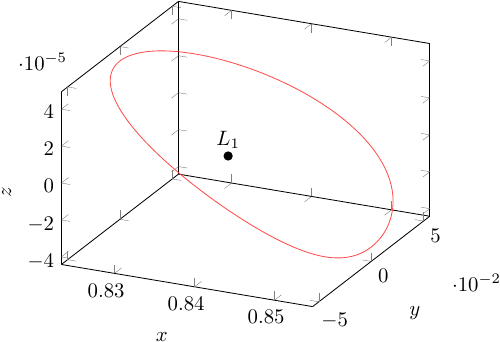}
  \includegraphics[page=20, scale=0.9]{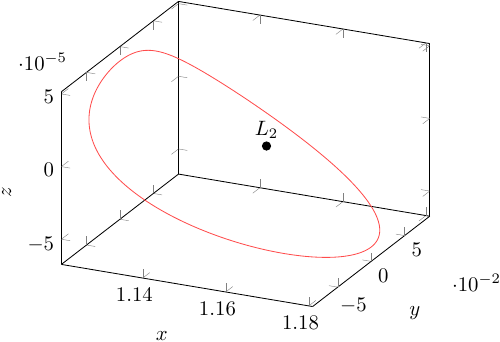}
  \caption{Halo orbits around $L_1$ and $L_2$ for energies
    \num{-1.5939} and \num{-1.5805} respectively. This corresponds to a unique
    excursion of \SI{16000}{km}.}
  \label{fig:Halo}
\end{figure}

\pardesc{Halo Orbits}
For the sake of
generality, we compute two Halo orbits with different energies. For the
Halo orbit around $L_{1}$ denoted by $\Halo_{1}$, we have chosen
$\energy(\Halo_{1})=\num{-1.5939}$.
For the Halo orbit around $L_{2}$ denoted by $\Halo_{2}$, we have
chosen
$\energy(\Halo_{2})=\num{-1.5805}$.
These two energy values correspond to a unique $z$-excursion of
$\SI{16000}{km}$.
The numerical computation of such orbits is done using the method
described in section~\ref{sec:Lyapunov} extended to the spatial case.
See figure~\ref{fig:Halo} for a plot of these two periodic orbits.

\subsubsection{Propagation of Manifolds and Choice of Trajectories}

\begin{figure}[ht]
  \centering
  \includegraphics[width=0.5\textwidth]{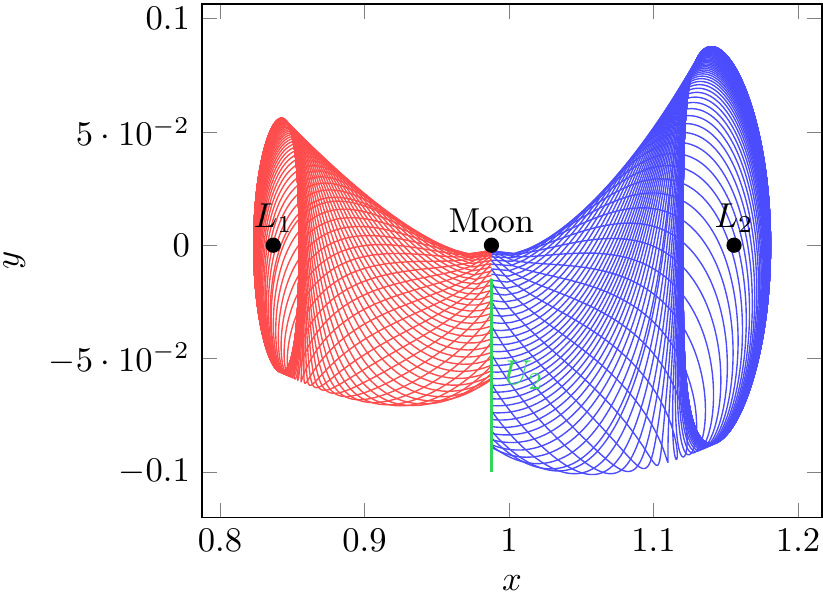}
  \includegraphics[page=21]{pdffigures-crop}
  \caption{On the left: the propagation of manifolds from Halo around
    $L_1$ and Halo around $L_2$. On the right: the two trajectories of
    these manifold minimizing the distance on the plane $U_2$. }
  \label{fig:HaloMani}
\end{figure}

Using the same parameter $\alpha$ as defined in~\eqref{eq:pertubLya}
for the Lyapunov to Lyapunov mission, \latin{i.e.} $\frac{1}{384402}$, we
compute the intersection with the plane $U_{2}$ (see
sec~\ref{sec:het}). One can see the result in
figure~\ref{fig:HaloMani}. We denote by $\Mani_{1}$ and $\Mani_{2}$
these two manifolds.

We compute the section of each manifold with $U_{2}$ and find the
closest pair of points (one from the manifold of $\Halo_{1}$
and one from the manifold of $\Halo_{2}$). This is done with a fine
discretization of 1000 points per Halo orbits. In that way, we get
two points for $x=1-\mu$,  denoted  respectively by
$\chi_{\Mani_{1}}^{U_{2}}$ and $\chi_{\Mani_{2}}^{U_{2}}$. The
distance in $\R^{6}$ is
  $\norm{\chi_{\Mani_{1}}^{U_{2}}- \chi_{\Mani_{2}}^{U_{2}}}{2}=\num{0.098644604436}$.
The two corresponding
trajectories are plotted in figure~\ref{fig:HaloMani}. Let
$t_{\Mani_{1}}$  and $t_{\Mani_{2}}$ denote  the
two times of propagation for the two free trajectories themselves
denoted by $A_{\Mani_{1}}$ and
$A_{\Mani_{2}}$ (see figure~\ref{fig:HaloMani}).

\subsubsection{Three Short Transfers}

Following our method, we compute three short transfers in order to
initialize a multiple indirect shooting method and get the optimal
trajectory.

\begin{figure}[ht]
  \centering
  \includegraphics[page=23]{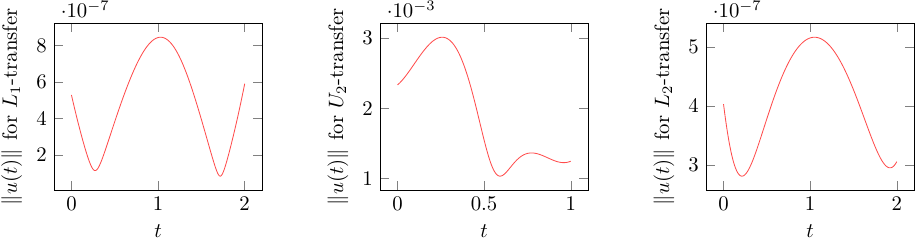}
  \caption{Norm of the control for the three controlled parts of the admissible trajectory.}
  \label{fig:localcontrol}
\end{figure}

\pardesc{From $\Halo_{1}$ to $A_{\Mani_{1}}$.} Once again, we follow the
 method described in section~\ref{sec:aroundL1}, we construct two fixed
points, one on $\Halo_{1}$, the other on $A_{\Mani_{1}}$. To do that,
we consider the two closest points on $\Halo_{1}$ and
$A_{\Mani_{1}}$.
We choose two time parameters: $t_{\Halo_{1}}$ to propagate backward
the point on
$\Halo_{1}$ and $t_{A_{\Mani_{1}}}^{L_{1}}$ to propagate forward on
$A_{\Mani_{1}}$. Here, we pick:
$t_{\Halo_{1}}=t_{A_{\Mani_{1}}}^{L_{1}}=\num{1.0}$.
We are now ready to build the first optimal control
problem as defined in~\eqref{eq:OCPL1}. Using continuation on the final
state solves this problem in an easy
and fast ($\SI{4.1}{s}$) manner.  The norm of the control is plotted in
figure~\ref{fig:localcontrol}. Let us denote by
$t_{0} =t_{\Halo_{1}}+t_{A_{\Mani_{1}}}^{L_{1}}$,
the transfer time and by $X_{0}^{*}=(\chi_{0}^{*},m_{0}^{*})$ and
$X_{1}=(\chi_{1},m_{1})$ the terminal
points of this transfer.  The resulting final mass  is
$m_{1}=\SI{1499.9967439278}{kg}$.

\pardesc{From $A_{\Mani_{1}}$ to $A_{\Mani_{2}}$.}
In this part, we apply our method to compute the transfer from
trajectory $A_{\Mani_{1}}$  to  trajectory
$A_{\Mani_{2}}$. There is a rather large  gap to
resorb. We already have the two points that we will perturb backward and
forward: $\chi_{\Mani_{1}}^{U_{2}}$ and $\chi_{\Mani_{2}}^{U_{2}}$. We
choose the two corresponding times
$t_{\Mani_{1}}^{U_{2}}=t_{\Mani_{2}}^{U_{2}}=\num{0.5}$,
and define the transfer time as\footnote{We keep the index $1$ for the
  remaining time on the free part, \latin{i.e.}, the remaining part of
  the unstable manifold trajectory.}
$t_{2} = t_{\Mani_1}^{U_{2}}+t_{\Mani_2}^{U_{2}}$.
After the first transfer from $\Halo_{1}$ to $A_{\Mani_{1}}$, we
follow a free trajectory on the manifold, so we choose the initial
mass of the
transfer from $A_{\Mani_{1}}$ to $A_{\Mani_{2}}$ as the final mass of
the previous part, that is to say $m_{2} =
m_{1}=\SI{1499.9967439278}{kg}$.

Once again, the continuation on the final state allows for a fast
convergence to obtain the solution of this problem. Indeed we obtained the
solution in $\SI{4.4}{s}$. In this problem, we denote by
$X_{2}=(\chi_{2},m_{2})$ the initial point on $A_{\Mani_{1}}$ and by
$X_{3}=(\chi_{3},m_{3})$ the final point on $A_{\Mani_{2}}$.
The final mass we get is $m_{3}=\SI{1493.3184622015}{kg}$ and the norm
of the control is plotted in
figure~\ref{fig:localcontrol}.

\pardesc{From $A_{\Mani_{2}}$ to $\Halo_{2}$.} We consider in here the
last short transfer from $A_{\Mani_{2}}$  to $\Halo_{2}$. Like the transfer from
$\Halo_{1}$ to $A_{\Mani_{1}}$, we pick the two closest points on
$A_{\Mani_{2}}$ and $\Halo_{2}$ and we perturb them with two time
parameters denoted by $t_{\Mani_{2}}^{L_{2}}$ and $t_{\Halo_{2}}$. We
define then the transfer\footnote{Once again we keep the index $3$
  for the free part between transfer between manifolds and transfer to
  $\Halo_{2}$.}
time $t_{4}=t_{\Mani_{2}}^{L_{2}}+t_{\Halo_{2}}=(1.0+1.0)$
to go from $X_{4}=(\chi_{4},m_{4})$ to the final points
$\chi_{5}^{*}$. As before, $m_{4}=m_{3}=\SI{1493.3184622015}{kg}$
because after the transfer around $U_{2}$, we follow a free trajectory
on the stable manifold.

Once again, the resolution is easy and fast thanks to the continuation
method: \SI{4.14}{s}. The final mass we obtain is $m_{5}
=\SI{1493.3156736966}{kg}$. The norm of the control is plotted in
figure~\ref{fig:localcontrol}.

\begin{figure}[ht]
  \centering
  \includegraphics[page=22]{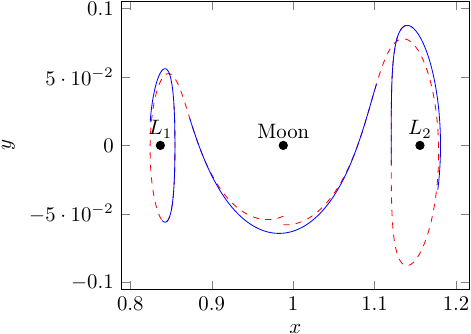}
  \caption{Admissible trajectory in five parts. Dashed part are manifold trajectories, \latin{i.e.} free parts.}
  \label{fig:admissibleHalo}
\end{figure}

\pardesc{Admissible Trajectory in Five Parts.} To summarize, we have
constructed an admissible trajectory going from $\Halo_{1}$ to
$\Halo_{2}$ with three controlled parts and two free parts. This
trajectory is plotted in the figure~\ref{fig:admissibleHalo} and the
control for the three controlled parts is in
figure~\ref{fig:localcontrol}. The local transfers are computed with a
maximal thrust equal to \SI{180}{N}, indeed this helps the
convergence of local transfers (but we do not reach the targeted
maximal thrust of \SI{0.3}{N}), and the multiple shooting
for the reason described in section~\ref{sec:aroundL1}.

\subsubsection{Multiple Shooting for the Total Transfer}
We consider here the following problem
\begin{equation}
  \label{eq:OCPTOTHalo}
  \Prob{\tot}^{\Halo}\left\{
    \begin{array}{l}
      \displaystyle \cost_{\tot}=\min \int_{0}^{t_{\tot}}\norm{u}{2}\D t,\\
      \displaystyle\dot x = F_{0}(x)+\frac{\epsilon}{m}\sum_{i=1}^{3}u_{i}F_{i}(x),\\
      \displaystyle\dot m = -\beta_{*}\epsilon\norm{u}{},\\
      \norm{u}{}\leq 1,\\
      x(0)=\chi_{0}^{*}\in\Halo_{1},\;m(0)=m_{0}^{*},\\[3pt]
      x(t_{\tot})=\chi_{5}^{*}\in\Halo_{2}.
    \end{array}
  \right.
\end{equation}
Transfer time $t_{\tot}$ is defined as
$t_{\tot}=t_{0}+t_{1}+t_{2}+t_{3}+t_{4},$
where $t_{0}$, $t_{2}$ and $t_{4}$ are the times previously introduced for
the three short transfers. Time $t_{1}$ is  the duration of
trajectory $A_{\Mani_{1}}$ in the unstable manifold from $\Halo_{1}$
from which we remove the two times we used to perturb points for the
local transfers around $L_{1}$ and $U_{2}$. This gives us
$t_{1}=t_{\Mani_{1}}-t_{A_{\Mani_{1}}}^{L_{1}}-t_{\Mani_{1}}^{U_{2}}$.
And we defined $t_{3}$ in a similar way as
$t_{3}=t_{\Mani_{2}}-t_{A_{\Mani_{2}}}^{L_{2}}-t_{\Mani_{2}}^{U_{2}}$.
As in section~\ref{sec:multipleshooting}, we introduced the shooting
function with four nodes (for the Lyapunov to Lyapunov mission, we
had two nodes). Because, we are considering a spatial mission, we get
the following shooting function
\begin{equation}
  \label{eq:multipleshootingHalo}
  \Ftir_{\multi}^{\Halo}(Z)=
  \begin{pmatrix}
    \flot^{\Ext}_{1,\dots,7}(\chi_{0}^{*},m_{0}^{*},P_{0})-X_{1}\\
    \flot^{\Ext}_{8,\dots,14}(\chi_{0}^{*},m_{0}^{*},P_{0})-P_{1}\\
    \flot^{\Ext}_{1,\dots,7}(X_{1},P_{1})-X_{2}\\
    \flot^{\Ext}_{8,\dots,14}(X_{1},P_{1})-P_{2}\\
    \flot^{\Ext}_{1,\dots,7}(X_{2},P_{2})-X_{3}\\
    \flot^{\Ext}_{8,\dots,14}(X_{2},P_{2})-P_{3}\\
    \flot^{\Ext}_{1,\dots,7}(X_{3},P_{3})-X_{4}\\
    \flot^{\Ext}_{8,\dots,14}(X_{3},P_{3})-P_{4}\\
    \flot^{\Ext}_{1,\dots,6}(X_{4},P_{4})-\chi_{5}^{*}\\
    \flot^{\Ext}_{14}(X_{4},P_{4})\\
  \end{pmatrix}\in\R^{4\times 14+7},
\end{equation}
where the vector $Z$ is defined as
\[Z=(\underbrace{p_{0},p_{m}^{0}}_{P_{0}},\underbrace{\chi_{1},m_{1}}_{X_{1}},
\underbrace{p_{1},p_{m}^{1}}_{P_{1}},\underbrace{\chi_{2},m_{2}}_{X_{2}},\underbrace{p_{2},p_{m}^{2}}_{P_{2}},\underbrace{\chi_{3},m_{3}}_{X_{3}},\underbrace{p_{3},p_{m}^{3}}_{P_{3}},\underbrace{\chi_{4},m_{4}}_{X_{4}},\underbrace{p_{4},p_{m}^{4}}_{P_{4}})\in\R^{63}.\]
We initialize the shooting method with the values that  we get from
the local transfers
and with a zero adjoint vector for the free parts. The shooting
converges easily.

\begin{figure}[ht]
  \centering
  \includegraphics[page=24]{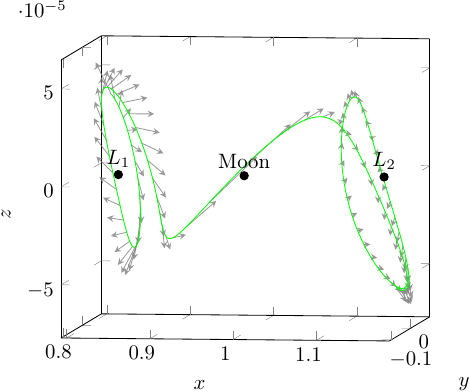}
  \includegraphics[page=25]{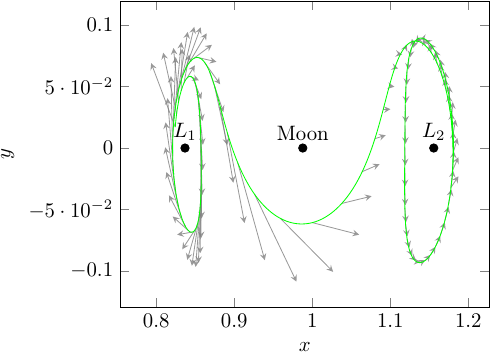}
  \caption{Optimal trajectory for the Halo to Halo transfer. On the
    left, a 3 dimensional view. On the right, a view in the
    $(x,y)$-plane. The control is represented by arrows.}
  \label{fig:HaloFinalTraj}
\end{figure}

\begin{figure}[ht]
  \centering
  \includegraphics[page=26]{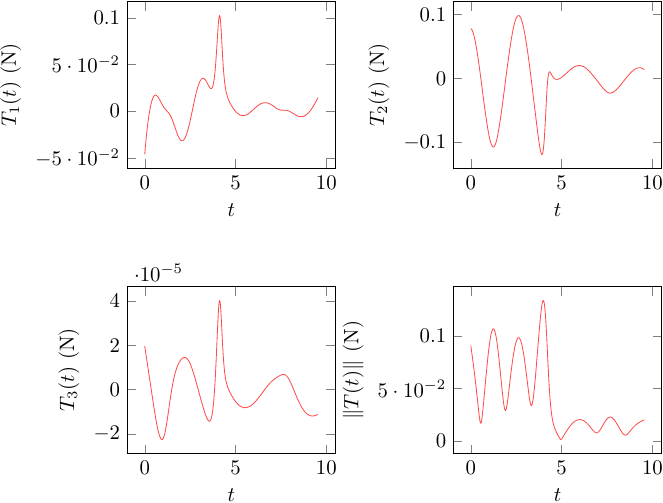}
  \caption{Optimal control for the Halo to Halo Mission $T(\cdot)\in\R^3$ in Newton.}
  \label{fig:HaloFinalControl}
\end{figure}

As for the two previous missions, we decrease the maximal authorized
thrust
 by  continuation. We then optimize the terminal
points $\chi_{0}^{*}$ and $\chi_{5}^{*}$ to satisfy the
transversality conditions. We manage to get the result in
\SI{4.16}{min}. The final trajectory is plotted in
figure~\ref{fig:HaloFinalTraj} and the corresponding control in
figure~\ref{fig:HaloFinalControl}. The result cost is summarized in
table~\ref{tab:HaloNumerical} as well as the numerical values of the
parameters. Because we are not comparing this mission with other
published results, we just write the physical cost $\cost_{\tot}^{3}$
in the international
system of units (see~\eqref{eq:costs}). Note that we do not get the
turnpike properties, indeed, in this case, there is no
``steady-state'' trajectory asymptotically connecting the two periodic
orbits.

\begin{table}[ht]
  \centering

  \begin{tabular}{ccc}
    \toprule
    Initial Mass &  Transfer time & $\Tmax$\\
    \midrule
    \SI{1500}{kg} & \num{9.5436454462828} or
    \SI{41.50}{days} & \SI{0.3}{N}\\
    \bottomrule
  \end{tabular}
  \bigskip

  \begin{tabular}{ccccc}
    \toprule
    & $\cost_{\tot}^{3}$  & Mass of fuel\\
    \midrule
    Halo to Halo Problem& \num{0.00461912647735513}& \SI{7.41587259099992}{kg} \\[5pt]
    \bottomrule
  \end{tabular}
  \bigskip

  \caption{Numerical results for the final trajectory of the Halo to
    Halo Mission.}
  \label{tab:HaloNumerical}
\end{table}

\section{Conclusion}
To design different spacecraft missions between periodic orbits around
Lagrange
points, we have used natural (uncontrolled) trajectories computed
thanks to the invariant
manifolds of the periodic orbits. We have connected resulting arcs
with short transfers using the PMP and indirect methods.
Doing that, we have designed admissible
trajectories performing the mission with controlled and uncontrolled
parts. The resulting admissible trajectories have been used to initialize
an indirect multiple shooting method in which we released the constraints
to join
uncontrolled parts, \latin{i.e.}, to force the spacecraft to follow the
natural drift. We have finally obtained a trajectory satisfying
the
first order necessary conditions for optimality given by the PMP.

In order to improve the robustness of our indirect approach, we have
designed and implemented appropriate
continuations on the final state and on the thrust. Thanks to this,
the execution of the overall computation is run within short time (of
order of a few minutes), and results have the excellent accuracy of
the underlying Newton method.

One can note that, when there is an heteroclinic orbit between the two
terminal periodic orbits, the optimal enjoys a turnpike property.
Proving the turnpike feature for such control-affine systems with
drift is an open issue, which may deserve consideration because it
gives an approach to successfully initialize a variant of the shooting
method in a simple and efficient way (see~\cite{trelatzuazua}).
Finally, as already mentioned, we have considered the
$L^{2}$-minimization of the cost, leaving the computation of the $L^{1}$-minimization solution with a
\emph{bang-bang} control as an open issue for farther studies.


\end{document}